%% file: helm_pv_paper_1.tex
\documentclass[journal, final, twoside, letterpaper]{IEEEtran}

\usepackage{colordvi}
\usepackage[dvipsnames]{xcolor}

\usepackage[cmex10]{amsmath}
\usepackage[T1]{fontenc}
\usepackage{letterspace}
\usepackage{amssymb}
\usepackage{amsfonts}
\usepackage{mathtools}
\usepackage{booktabs}
\usepackage{psfrag}
\usepackage{multirow}
\usepackage{array}
\usepackage{multirow} 
\usepackage{fixltx2e}
\usepackage{cite}
\usepackage{graphicx}
\usepackage{commath}
\usepackage[font=footnotesize,caption=false]{subfig}
\usepackage{pgfplots}
\usepackage{tikz}
\usetikzlibrary{shapes,fit,arrows,decorations.markings,matrix}
\usepackage[mediumspace,mediumqspace,squaren,binary]{SIunits}
\usepackage{threeparttable}
\pgfplotsset{
	compat=1.3,
	/pgfplots/ylabel absolute/.style={%
		/pgfplots/every axis y label/.style={at={(0,0.5)},xshift=-15pt,rotate=90},
		/pgfplots/every y tick scale label/.style={at={(0,1)},above right,inner sep=0pt,yshift=0.3em}
	}
}%

\newcommand{\mc}{\mathcal}

\newcommand{\sumoN}{\sum_{k\in\mc{B}}}

\newtheorem{assumption}{Assumption}

\title{Alternative PV Bus Modelling with the Holomorphic Embedding Load Flow Method}%
\author{I.~Wallace, D.~Roberts,
  A.~Grothey, and K.~I.~M.~McKinnon --- \today%
\thanks{I.~Wallace, A.~Grothey, and K.~I.~M.~McKinnon 
are with the
  School of Mathematics, University of Edinburgh, James Clerk Maxwell
  Building, Edinburgh EH9~3JZ, UK (e-mail:
  \texttt{i.p.wallace@sms.ed.ac.uk},
  \texttt{a.grothey@ed.ac.uk}, \texttt{k.mckinnon@ed.ac.uk}).}}%
\begin{document}
\maketitle

\begin{abstract}
  The Holomorphic Embedding Load Flow Method (HELM) has been suggested
  as an alternative approach to solve load flow problems. However, the
  current literature does not provide any HELM models that can
  accurately handle general power networks containing PV and PQ buses
  of realistic sizes.  The original HELM paper dealt only with PQ
  buses, while a second paper showed how to include PV buses but
  suffered from serious accuracy problems.  This paper fills this gap
  by providing several models capable of solving general networks,
  with computational results for the standard IEEE test cases provided
  for comparison. In addition this paper also presents a new
  derivation of the theory behind the method and investigates some of
  the claims made in the original HELM paper.

\end{abstract}


\section{Introduction}

\IEEEPARstart{L}{oad} flow problems consist of solving a set of
equations --- the bus-power-equilibrium equations (BPEE) --- based on
the physical laws surrounding a power network so as to obtain a
feasible load flow for the network. The physical laws result in the
BPEE being nonlinear, and usually having multiple solutions.  In
addition, many of these solutions would typically correspond to 
{\em unstable} operating points of a real network, so termed {\em low
voltage solutions}, and so solving the BPEE for the physical solution is not an
 easy calculation.

 Most current methods for solving the BPEE of a network use an
 iterative approach. These methods all have the same issue --- that
 for specific problems they may not converge to a solution or converge
 to the ``wrong'' ({\em i.e.,} unstable or low voltage) solution. A 2012
 paper \cite{helm} and a subsequent patent \cite{helmpatent} by Trias
 proposed an alternative to these iterative methods, HELM, which
 claims to address this issue. The idea was to treat the voltages as
 holomorphic functions of a complex parameter $z$ that scales the
 demands, and to use the simple-to-calculate solution when $z=0$ to
 determine the desired solution when $z=1$ making use of the powerful
 theory of holomorphic functions and analytic continuations. There
 remained some important gaps in the method as described by Trias
 however, as it did not provide any detail on how to handle PV buses,
 and so the method was applicable only on networks with exclusively PQ
 buses (plus one slack bus). As well, Trias' papers contain
 claims about the theoretical properties of the method, namely
 regarding the method's ability to always find the ``correct'' solution,
 for which only sketchy theoretical justification was given.
 Subramanian {\em et al} \cite{subramanian} suggest an extension of
 the HELM approach that can be applied to networks that include PV buses,
 however this approach suffers from numerical problems.

 The contribution of this paper is to provide alternative approaches
 to model PQ buses within the HELM framework that do not suffer from
 the above numerical problems. Further, we give an investigation of the
 theoretical underpinnings of the HELM method.

 The organisation of this paper is as follows. In the following
 section, we recap the HELM method for PQ bus networks as proposed in
 \cite{helm,helmpatent}. In Section~\ref{sec:theory}, we give a
 derivation of the HELM approach that uses the holomorphic implicit
 function theorem. This enables us to gain further insight into the 
 theoretical properties of the method. In Section~\ref{sec:subramanian}, 
 we show how the problem changes when considering PV buses, present the 
 current literature and its shortcomings, and then in 
 Section~\ref{sec:PVmodels} provide our own models which are capable 
 of handling general networks. Section~\ref{sec:results} presents 
 computational results for the test networks.
 
\section*{Nomenclature}

\subsection*{Sets}
\begin{description}
\item[$\mc{B}$] All Buses.
\item[$\mc{B}_{PQ}$] PQ Buses.
\item[$\mc{B}_{PV}$] PV Buses.
\end{description}
\subsection*{Parameters and Variables}
\begin{description}
\item[$Y_{ik}$] ($i,k$) element of bus admittance matrix.
\item[$V_i$] Complex voltage at bus $i$.
\item[$S_i$] Complex power injection at bus $i$ ($S_i=P_i+jQ_i$).
\item[$P_i$] Real power injection at bus $i$.
\item[$Q_i$] Reactive power injection at bus $i$.
\item[$M_i$] Prescribed voltage magnitude at PV bus $i$.
\item[$z$] Complex variable used for holomorphic embedding.
\end{description}
\subsection*{Notation}
\begin{description}
\item[$\delta_{i,j}$] $1$ if $i\!=\!j$, $0$ otherwise.
\end{description}

\section{HELM PQ Model}
 \label{sec:helm}
 The HELM method was first introduced by Trias in \cite{helm}. He
 begins by considering the BPEE for a PQ-bus network in the general
 form:
 \begin{align}
 \sum_{k\in \mc{B}}Y_{ik}V_k=\frac{S_i^*}{V_i^*},\quad i \in \mc{B}_{PQ}\label{eq:bpee}
 \end{align}
 where without loss of generality we will set bus 0 to be the slack
 bus, so that $\mathcal{B} = \mathcal{B}_{PQ}\cup\{0\}$. The voltage at the slack bus is known to be $V_0=1$. The remaining $V_i$ are
 the unknown complex variables. Trias proceeds by setting up a homotopy where the demands are scaled by a {\em complex} parameter $z$ and the resulting bus voltages are treated as functions of this complex parameter
 \begin{align}
 \sum_{k\in \mc{B}}Y_{ik}V_k(z)=\frac{zS_i^*}{V_i^*(z^*)},\quad i \in \mc{B}_{PQ}\label{eq:bpeez}
 \end{align}
 here  $z\!=\!1$ corresponds to the solution of the BPEE while 
 at $z\!=\!0$ a solution can be easily computed.
Note that the solution for \eqref{eq:bpeez} will in general not be unique for a given $z$, however for $z=0$ a unique solution exists under the condition that $V_i(0)\ne 0\,\forall i$. 

The main claim by Trias is that the voltages $V_i(z)$ implicitly
defined by \eqref{eq:bpeez} are holomorphic functions at $z=0$ and can
be analytically continued to obtained the ``correct'' solution to the
BPEE at $z=1$.\footnote{To be precise: The unique seed, or reference, solution at $z=0$ under the
  condition $V_i(0)\ne 0$ can be continued into a neighbourhood by the
  implicit function theorem. This continuation is holomorphic at
  $z=0$.}  This is not obvious, due to the use of the complex conjugate in system \eqref{eq:bpeez} which is {\em not} a holomorphic function. 

Trias circumvents this difficulty by embedding \eqref{eq:bpeez} in a
larger holomorphic system, namely

 \begin{equation}\begin{aligned}
 \sum_{k\in \mc{B}}Y_{ik}V_k(z)=\frac{zS_i^*}{\overline{V_i}(z)},\quad i \in \mc{B}_{PQ}\\
 \sum_{k\in \mc{B}}Y^*_{ik}\overline{V_k}(z)=\frac{zS_i}{V_i(z)},\quad i \in \mc{B}_{PQ}
 \end{aligned}\label{eq:helm-first}\end{equation}
where $\overline{V_i}(z)$ are additional complex variables formally
independent of $V_i(z)$. It is easy to check that these equations are
indeed holomorphic as functions of the independent complex variables
$z, V_i, \overline{V_i}$ for example by checking the Wirtinger
derivatives or the Cauchy-Riemann equations.
System \eqref{eq:helm-first} is a set of polynomial equations
(after multiplying through with the denominator in each case) and
Trias uses the theory of resultants and Gr\"obner bases to deduce that
all $V_i$ and $\overline{V_i}$ are holomorphic functions everywhere
except for a finite set of singularities - all of them branch points - which will not include 0. 

If the additional constraint
\begin{equation}
\overline{V_i}(z) = (V_i(z^*))^*, \quad i\in \mc{B},\label{eq:helm-refcond}
\end{equation}
which Trias calls the {\em reflecting condition}, holds, system
\eqref{eq:helm-first} reduces to \eqref{eq:bpeez}.  
Trias makes use of
the system \eqref{eq:helm-first} only to establish that there exist
holomorphic solution functions $V_i(z), \overline{V_i}(z)$ and then
argues that since we are only interested in those solutions that
satisfy the reflecting condition it can be used to eliminate the
$\overline{V_i}$. In the remainder of his presentation Trias uses \eqref{eq:bpeez} exclusively.

In Section \ref{sec:theory} we will show that the reflecting condition
\eqref{eq:helm-refcond} is satisfied automatically under the condition
that $V_i(0)\ne 0\,\forall i$.

Given that the voltages are holomorphic functions in a neighbourhood
of $z\!=\!0$, they, and their reciprocals, can be written as power series expandable about
$z\!=\!0$:
\begin{align}
V_i(z)&=\sum_{n=0}^{\infty} V_i[n]z^n, \quad i\in \mc{B}\\
\frac{1}{V_i(z)}=W_i(z)&=\sum_{n=0}^{\infty} W_i[n]z^n, \quad i\in \mc{B}. \label{eq:helm-Wdef}
\end{align}
Substituting into \eqref{eq:bpeez} one obtains
\begin{align}
 \sum_{k\in\mc{B}}Y_{ik}\sum_{n=0}^{\infty} V_k[n]z^n= zS_i^*\sum_{n=0}^{\infty} W_i^*[n]z^n\label{eq:helm-powseries}
 \end{align}
From \eqref{eq:helm-powseries}, it is now possible to determine the coefficients of the power series up to any desired level. The process is begun by solving for $z\!=\!0$, which (under the condition that $V_i(0)\ne 0$) yields the set of linear equations
\begin{align}
\sum_{k\in \mc{B}}Y_{ik}V_k[0]=0, \quad i \in \mc{B}_{PQ}.\label{eq:helm-v0}
\end{align}
Note that the sum on the left includes the slack bus $k=0$ for which
$V_0(z)$ is set to $1$ for all $z$. 
At this point we need to impose
\begin{assumption}
The reduced bus admittance matrix $Y'$ obtained from $Y$ by removing the row and column corresponding to the slack bus is non-singular. 
\end{assumption}
This is a standard assumption and will hold for any sensible power system. In particular, in the absence of shunts and phase shifters the assumption is equivalent to requiring the system to be connected. 

Under Assumption~1 system (\ref{eq:helm-v0}) has a unique solution.
$W_i[0]$ can then be computed using \eqref{eq:helm-Wdef}
\begin{equation}
W_i[0]=\frac{1}{V_i[0]}
\end{equation}
Having obtained the initial values for $V$ and $W$, an iterative process can be used to determine the remaining values in the power series up to any desired order of $n$ by equating the coefficients of $z^n$ in \eqref{eq:helm-powseries}, which yields
\begin{align}
\sum_{k\in \mc{B}}Y_{ik}V_k[n]= S_i^*W^*_i[n\!\!-\!\!1], \quad i \in \mc{B}_{PQ}\quad n\geq 1\label{eq:helm-vcalc}
\end{align}
where $W_i[n\!-\!1]$ is calculated using the coefficients of lower orders
\begin{align}
W_i[n-1]=-\frac{\sum\limits_{m=0}^{n-2}V_i[n\!-\!m\!-\!1]W_i[m]}{V_i[0]}\label{eq:helm-wcalc}
\end{align}
In \eqref{eq:helm-v0} and \eqref{eq:helm-vcalc} the coefficient matrix of the system of equations is constant, and so factorisation of this matrix needs only to be done once and can be used for all iterations. 

Having obtained the power series for the voltages up to some desired level, it is now possible to compute the voltages for each bus. However, a direct summation of the power series for $z\!=\!1$ is insufficient, as the radius of convergence of the power series is typically less than 1. Instead, analytic continuation using Pad\'e approximants \cite{pade} is used.
Pad\'e approximants are a particular type of rational approximations to power series known to have good convergence properties. The $L,M$ Pad\'e approximant is denoted by $[L/M]=P_L(x)/Q_M(x)$, where $P_L(x)$ and $Q_M(x)$ are polynomials of degree less than or equal to $L$ and $M$ respectively.

In \cite{helm}, Trias explains how Stahl's extremal domain theorem and Stahl's Pad\'e convergence theorem provide proof that Pad\'e approximants give the maximal analytical continuation. 
That is, if there is a steady-state solution to the problem, then the diagonal Pad\'e approximants will converge to this answer, while if there is no steady-state solution (voltage collapse), then the Pad\'e approximants will not converge.
In fact Stahl's
Theorems \cite[Ch 6]{pade} asserts that the diagonal Pad\'e
approximants ({\em i.e.,} $L=M$) converge to the analytic continuation
with the extremal domain of the approximated function. Here extremal
domain is understood as the one having a minimal exemption set (in the
shape of branch cuts) measured in capacity. The implication is that diagonal Pad\'e approximants (of high enough order) will yield values $V_i(1)$ provided that $z=1$ is not on a branch cut. Such a branch cut is indicated by a line of poles of the approximant.

A slight variant to \eqref{eq:helm-first} is given in \cite{helmpatent}, where an additional term is added to the holomorphic embedding
\begin{equation}\begin{aligned}
 \sum_{k\in \mc{B}}Y_{ik}V_k(z)-(1\!\!-\!\!z)y_i=\frac{zS_i^*}{\overline{V_i}(z)},\quad i \in \mc{B}_{PQ}\\
 \sum_{k\in \mc{B}}Y^*_{ik}\overline{V_k}(z)-(1\!\!-\!\!z)y_i^*=\frac{zS_i}{V_i(z)},\quad i \in \mc{B}_{PQ}
 \end{aligned}\label{eq:helm-alternative}\end{equation}
where
\begin{align}
y_i=\sum_{k\in\mc{B}}Y_{ik}
\end{align}

In this alternative model, at $z\!=\!0$, $V_k(0)=1$ is an obvious solution, as this causes the two terms in the left-hand side of the equations to cancel. This also means that $V_k[0]=1$, eliminating the need for step \eqref{eq:helm-v0}. At $z\!=\!1$ the $(1\!-\!z)$ term disappears, leaving the equations the equivalent of \eqref{eq:helm-first}. This alternative model will be the basis of some of the models which incorporate $PV$ buses shown later on.

\section{Theory}
\label{sec:theory}
In this section we will present additional theory for the HELM model. First, we will provide a separate proof that the $V_i$ and $\overline{V_i}$ in \eqref{eq:helm-first} are holomorphic using the Complex Implicit Function Theorem (CIFT) \cite{ift}, which we will then extend in Section \ref{sec:PVmodels} to include models with PV buses. Next we will show that the reflecting condition \eqref{eq:helm-refcond} is implied by the formulation and thus need not be assumed.

\subsection{$V$ and $\overline{V}$ Holomorphic}\label{sec:HELM-holo}
 We begin by defining the functions
\begin{equation}\begin{aligned}
f_i(z,V,\overline{V}):=&\overline{V_i}\sumoN Y_{ik}V_k - zS_i^*,\quad &i\in \mc{B}_{PQ}\\
f_{N+i}(z,V,\overline{V}):=&V_i\sumoN Y^*_{ik}\overline{V_k} - zS_i, \quad &i\in\mc{B}_{PQ}
\end{aligned}\label{eq:ift}\end{equation}
where $N$ is the number of non-slack buses in the network.

The CIFT states that if there exists a seed solution $v, \overline{v}$ with $f(0,v,\overline{v})=0$ and $J$ is non-singular at $(0,v,\overline{v})$, where
\begin{equation}
J_{ij}=\frac{\partial f_i}{\partial U_j}, \quad i,j = 1\dots 2N,
\label{eq:CIFTJ}\end{equation}
and $U~:=~\{V_1,\dots , V_N,\overline{V}_1,\dots , \overline{V}_N\}$, then there exist holomorphic functions $V_i(z)$ and $\overline{V_i}(z)$ of $z$ that satisfy (\ref{eq:ift}) in a neighbourhood of $z\!=\!0$. 

In this setup, $f$ is clearly a holomorphic mapping, and the values of $v$ in the seed solution $f(0,v,\overline{v})=0$ are the solution to \eqref{eq:helm-v0}. 

Using \eqref{eq:ift}, the values of $J$ can be computed as follows
\begin{equation}\begin{aligned}
\frac{\partial f_i}{\partial V_j}(0,v,\overline{v})=\overline{v_i}Y_{ij},\quad i,j=1,\dots,N\\
\frac{\partial f_i}{\partial \overline{V_j}}(0,v,\overline{v})=0,\quad i,j=1,\dots,N\\
\frac{\partial f_{N+i}}{\partial V_j}(0,v,\overline{v})=0,\quad i,j=1,\dots,N\\
\frac{\partial f_{N+i}}{\partial \overline{V_j}}(0,v,\overline{v})=v_iY^*_{ij},\quad i,j=1,\dots,N
\end{aligned}\end{equation}
and so $J$ can be rewritten as
\begin{equation}
J=\begin{pmatrix}\overline{v'}Y'&0\\0&v'(Y')^*\end{pmatrix}
\end{equation}
where $v'=\textrm{diag}(v_1,\dots,v_N)$ and $Y'$ represents the admittance matrix without the slack bus row and column. 
Clearly $J$ is non-singular iff $Y'$ is non-singular which is guaranteed by Assumption~1 and therefore $V(z)$ and $\overline{V}(z)$ are holomorphic functions of $z$.

\subsection{Reflecting Condition Redundancy}\label{sec:refcond1}
Now we will show that, in a neighbourhood of $z=0$, the reflecting condition \eqref{eq:helm-refcond} is redundant. In particular we will show that any solution to (\ref{eq:helm-first})---which will automatically satisfy $V_i(0)\ne 0\,\forall i$---must satisfy the reflecting condition.
If we do not use the reflecting condition, then instead of \eqref{eq:helm-powseries} we obtain the following set of equations
\begin{equation}\begin{aligned}
 \sum_{k\in\mc{B}}Y_{ik}\sum_{n=0}^{\infty} V_k[n]z^n&= zS_i^*\sum_{n=0}^{\infty} \overline{W_i}[n]z^n\\
 \sum_{k\in\mc{B}}Y_{ik}^*\sum_{n=0}^{\infty} \overline{V_k}[n]z^n&= zS_i\sum_{n=0}^{\infty} W_i[n]z^n
 \end{aligned}\label{eq:refcondred}\end{equation}
where again $W_i=1/V_i$ and likewise $\overline{W_i}=1/\overline{V_i}$. In comparing coefficients of $z^n$, we get:
\begin{equation}
\begin{bmatrix}Y'&0\\0&(Y')^*\end{bmatrix}\begin{bmatrix}V[n]\\ \overline{V}[n]\end{bmatrix} = \begin{bmatrix}r_{1,n-1} \\ r_{2,n-1}\end{bmatrix}
\label{eq:refcondred1}
\end{equation}
with 
\begin{equation}\begin{aligned}
r_{1,n-1}&=\left\{\begin{array}{lr}
 -Y_0&n=0\\
S_i^*\overline{W_i}[n-1]&n\geq 1 
 \end{array}\right.\\
r_{2,n-1}&=\left\{\begin{array}{lr}
 -Y_0^*&n=0\\
S_iW_i[n-1]&n\geq 1 
 \end{array}\right.\\
 \end{aligned}\label{eq:refcondrhs1}\end{equation}
where $Y_0$ is the slack bus column of the admittance matrix.

Now if we take the complex conjugate of the above system, we obtain:
\begin{equation}
\begin{bmatrix}(Y')^*&0\\0&Y'\end{bmatrix}\begin{bmatrix}V^*[n]\\ \overline{V}^*[n]\end{bmatrix} = \begin{bmatrix}r_{1,n-1}^* \\ r_{2,n-1}^*\end{bmatrix}
\end{equation}
and by rearranging
\begin{equation}
\begin{bmatrix}Y'&0\\0&(Y')^*\end{bmatrix}\begin{bmatrix}\overline{V}^*[n]\\ V[n]\end{bmatrix} = \begin{bmatrix}r_{2,n-1}^* \\ r_{1,n-1}^*\end{bmatrix}
\label{eq:refcondred2}\end{equation}

From \eqref{eq:refcondrhs1} it follows that when $n\!=\!0$, $r_{2,-1}=-Y_0^*=r_{1,-1}^*$. Thus the right-hand side in \eqref{eq:refcondred2} is the same as in \eqref{eq:refcondred1}. Under Assumption~1 $Y'$ is non-singular and therefore the solutions of the two systems must be the same: namely $V[0]=\overline{V}[0]^*$.

We now assume that $V[n]=\overline{V}^*[n]$ --- and hence $W[n] = \overline{W}[n]^*$ --- is true up to $n=k$ and check that for $n=k+1$

\begin{equation}\begin{aligned}
 r_{1,k}&= S^*\overline{W}[k] \\
 &=S^*W[k]^*&\textrm{from inductive hypothesis}\\
 &=r_{2,k}^*
 \end{aligned}\end{equation}
Once again we have that the right hand side in \eqref{eq:refcondred2} is identical to the one in \eqref{eq:refcondred1}, and so $V[n]=\overline{V}[n]^*$ for $n\geq0$. The reflecting condition holds for all coefficients of the power series of $V_i(z)$ and $\overline{V}_i(z)$ and therefore $\overline{V}_i(z)=V_i(z^*)^*$ holds for the functions themselves. Moreover, if the [L/M] Pad\'e Approximant for $V(z)$ is given by $p(z)/q(z)$, then it is straightforward to show that the [L/M] Pad\'e Approximant to $\overline{V}(z)$ is given by $(p(z^*))^*/(q(z^*))^*$ and so the reflecting condition holds also for the analytically continued function. From, for example, \cite{Ahlfors} we know that the analytically continued functions satisfy the polynomial HE equations \eqref{eq:ift}

\section{Introducing PV Buses}
\label{sec:subramanian}
The HELM method as described in \cite{helm} had one major deficiency --- it did not describe how to deal with networks that include PV buses.
When considering PV buses, as seen in Table \ref{tab:bustypes}, the unknowns in the BPEE are different. At a PQ bus the real and reactive power injections are known and the complex voltage is unknown, whereas at a PV bus the real power injection and the voltage magnitude are known and the reactive power injection and the voltage angle are unknown. As we are solving for different variables in the BPEE, it is necessary to rearrange the equations. Without loss of generality we will consider the systems to be ordered such that the PQ buses are grouped first and then the PV buses come afterwards.
\begin{table}[htbp]
\centering\footnotesize
\caption{Types of buses in a power system}
\begin{tabular}{ccccc}
\toprule
Type of& Voltage Mag & Voltage &Real Power&Reactive Power\\
Bus&(|V|)&Angle ($\delta$)&Injection (P)&Injection (Q)\\
\midrule
Slack (V$\delta$)&Given&Given&Unknown&Unknown\\
PQ Bus&Unknown&Unknown&Given&Given\\
PV Bus&Given&Unknown&Given&Unknown\\
\bottomrule
\end{tabular}
\label{tab:bustypes}
\end{table}

In \cite{subramanian}, Subramanian et al present an approach to deal with general networks that may include PV buses. For PQ buses, their holomorphic embedding of the BPEE uses the alternative set of equations \eqref{eq:helm-alternative}.

For PV buses, the voltage magnitude and real power at the bus are known, but not the reactive power. Thus the authors create equations that use only the real power at a bus. Adding a number's complex conjugate to itself results in eliminating the imaginary part of that number---using this strategy with equation \eqref{eq:bpee}, the authors suggest the following holomorphic embedding to replace (\ref{eq:helm-alternative}) for PV buses
\begin{equation}\begin{aligned}
M_i^2\sumoN Y_{ik}V_k(z)& = z2P_iV_i(z) + (1\!\!-\!\!z)M_i^2y_i - {}\\
&z\!\left(V_i(z)^2\sum_{k=0}^NY_{ik}^*\overline{V_k}(z)\right) , i\in \mc{B}_{PV}\\
\overline{V_i}(z)V_i(z) &= 1+z(M_i-1),\ \ \ i\in\mc{B}_{PV}
\end{aligned}\label{eq:mod1-syseqnPV}
\end{equation}
where $M_i$ is the target voltage magnitude for PV buses.
The second equation in the holomorphic embedding, is not explicitly shown in \cite{subramanian}, but is required to provide a path for the voltage magnitudes in the PV buses to start at 1 when $z\!=\!0$ and finish at $M_i$ when $z\!=\!1$.

Now the reflecting condition \eqref{eq:helm-refcond} substituted into \eqref{eq:helm-alternative} and \eqref{eq:mod1-syseqnPV} combine to form the holomorphic embedding of the entire network.

When $z\!=\!0$ a solution to the system is simply $V_i(0)=1$, $i\in\mc{B}$. The power series coefficients are then determined by the same process as in the original HELM method. However, for PV buses, the term $V_i(z)^2\sum_{k=0}^NY_{ik}^*\overline{V_k}(z)$ contains products of three power series. This results in double convolutions, which can have precision limitations and can lead to inaccuracies in the final results. This problem is discussed in~\cite{subramanian}.

When applied to the IEEE test cases, the Subramanian model in \cite{subramanian} had poor convergence in even the simplest 9-Bus case, which is confirmed by our results in Table \ref{tab:1padecase9}. Here, $Rs$ refers to the residual when the voltages are substituted into the original BPEE, and $\Delta$ refers to the difference between the model voltage results and the voltage results obtained through MatPower. For more complicated cases the model was unable to provide even approximately correct results, as shown for the 39-Bus case in Table \ref{tab:1padecase39}. 

\begin{table}
\caption{Results for IEEE case9 using Subramanian Model}
\label{tab:1padecase9}
\centering
\begin{tabular}{ccc}
\toprule
Pad\'e Order& Max $|Rs|$ & Max $|\Delta|$ \\
\midrule
$[5/5]$ & 1.4e-01 & 1.2e-01 \\
$[10/10]$ & 3.9e-01 & 7.4e-02 \\
$[15/15]$ & 2.2e-02 & 4.8e-03 \\
$[20/20]$ & 5.4e-03 & 7.9e-04 \\
$[25/25]$ & 4.8e-04 & 6.6e-05 \\
\bottomrule
\end{tabular}
\end{table}

\begin{table}
\centering
\caption{Results for IEEE case39 using Subramanian Model}
\label{tab:1padecase39}
\begin{tabular}{ccc}
\toprule
Pad\'e Order& Max $|Rs|$  & Max $|\Delta|$ \\
\midrule
$[5/5]$ & 2.3e+01 & 4.2e+00 \\
$[10/10]$ & 4.8e+00 & 1.1e+00 \\
$[15/15]$ & 1.2e+01 & 8.5e-01 \\
$[20/25]$ & 2.5e+01 & 2.1e+00 \\
$[25/25]$ & 1.1e+01 & 1.1e+00 \\
\bottomrule
\end{tabular}
\end{table}

Updating his previous work in \cite{subramanian}, Subramanian provides a revised model in his thesis \cite{subrthesis}. This revised model no longer has the double convolution issue and is better suited to solving larger models. The model also no longer involves a $(1-z)$ term, instead it splits the admittance matrix into two parts, creating a diagonal matrix for the shunt effects which allows the remaining transmission elements to have zero row sums. The shunt elements are moved to the right-hand side and are multiplied by the complex parameter $z$.\footnote{though for simplicity he does not model the shunt conductance term for PV buses in his paper} The holomorphic embedding becomes

\begin{equation}\begin{aligned}
&\!\!\!\!\sumoN Y_{\text{trans}_{ik}}V_k(z) = \frac{zS_i^*}{V_i^*(z^*)} -zY_{\text{shunt}_{ii}}V_i(z),\ i\in \mc{B}_{PQ}\\
&\!\!\!\!\left(\!\!V_i(z)\sumoN Y_{\text{trans}_{ik}}^*V_k^*(z^*)\!\!\right)\!\!+\!\!\left(\!\!V_i^*(z^*)\sumoN Y_{\text{trans}_{ik}}V_k(z)\!\!\right) \\
&\ \ = 2zP_i - 2\text{Re}\left\{Y_{\text{shunt}_{ii}}V_i(z)V_i^*(z^*)\right\},\ \ \ i\in\mc{B}_{PV}\\
&\!\!\!\!\overline{V_i}(z)V_i(z) = 1+z(M_i-1),\ \ \ i\in\mc{B}_{PV}
\end{aligned}\label{eq:subrthes-syseqnPQ}
\end{equation}

This has the similar effect of making the seed voltages at every bus equal to $1$.  This model is different from our general model, which we provide in the next section. By separating the shunt elements from the admittance matrix, it creates an additional term which is dependent on the square of the voltage. In contrast, our general model is capable of creating either a linearly voltage-dependent term or a voltage-independent term. The results of both of Subramanian's models are given for the IEEE test cases in Section \ref{sec:results}. 

\section{New PV Models}
\label{sec:PVmodels}
In this section we present a general parametrised model for mixed PQ/PV systems that does not have the drawback of having a double convolution. We will show that for models of this general form the $V$ and $\overline{V}$ are holomorphic functions in $z$ and that the reflecting condition \eqref{eq:helm-refcond} remains redundant. Finally we will present specific choices for the parameters of the general model and investigate the numerical behaviour of the resulting methods. 

\subsection{General Model}
We suggest the following general model: for PQ buses, the following equations are used in the place of (\ref{eq:helm-alternative}):
\begin{equation}\begin{aligned}
\overline{V_i}(z)  \left( \sumoN Y_{ik}V_k(z)+ (z\!\!-\!\!1)a_i \right) &+{}\\
(z\!\!-\!\!1)b_i &= zS_i^*\quad i\in \mc{B}_{PQ}\\
V_i(z)\left( \sumoN Y^*_{ik}\overline{V_k}(z)+ (z\!\!-\!\!1)a^*_i \right) &+{}\\
(z\!\!-\!\!1)b^*_i &= zS_i\quad i\in \mc{B}_{PQ}
\end{aligned}\label{eq:genPQ}\end{equation}

This is similar to the approach used in both \cite{helmpatent} and \cite{subramanian}, and is identical if $a_i = y_i$ and $b_i=0,\quad i\in \mc{B}_{PQ}$. Natural choices for $a_i$ and $b_i$ would be $y_i$ or $0$ and these terms would serve a similar purpose as the additional $(1-z)y_i$ term in (\ref{eq:helm-alternative}). $a_i$ is multiplied by the voltage-dependent term, while $b_i$ is voltage-independent. For PV buses we use the following equations in the place of \eqref{eq:mod1-syseqnPV}:
\begin{equation}\begin{aligned}
\overline{V_i}(\!z\!)\!\!&\left(\sumoN\!\! Y_{ik}V_k(\!z\!)\!+\!(z\!\!-\!\!1)a_i\!\!\right)\!+\! V_i(\!z\!)\!\!\left(\sumoN \!\!Y^*_{ik}\overline{V_k}(\!z\!)+{}\right.\\
&\left.\vphantom{\sumoN}(z\!\!-\!\!1)a^*_i\right)+(z\!\!-\!\!1)(b_i+b_i^*) = 2zP_i\quad i\in \mc{B}_{PV}
\end{aligned}\label{eq:genPV1}\end{equation}
\begin{align}
\overline{V_i}(z)V_i(z)=L_i^2(z)\quad i\in \mc{B}_{PV}\label{eq:genPV2}
\end{align}
where $L_i(z)$, which is assumed to be real-valued, describes how the target voltage magnitude changes with respect to $z$.
\\ 
Equation \eqref{eq:genPV1} uses the same approach to eliminate the unknown reactive power as shown in \cite{subramanian}, except there is one less voltage term in both the left and right hand sides, so our model does not give rise to any double convolutions and their possible accuracy problems. Equation \eqref{eq:genPV2} gives us the voltage magnitude at each PV bus, which may depend on the value of $z$. The form of $L_i(z)$ is unrestricted except that $L_i(1)$ must equal $M_i$.

Equating coefficients of $z^n,n\geq 1$ in \eqref{eq:genPQ} will yield the following equations for $n=1$
\begin{equation}\begin{aligned}
\overline{V_i}[0]\sumoN &Y_{ik}V_k[1] + \overline{V_i}[1]\left(\sumoN Y_{ik}V_k[0]-a_i\right)= \\
&S_i^*-\overline{V_i}[0]a_i-b_i, \qquad\qquad\qquad\quad i \in \mc{B}_{PQ}\\
V_i[0]\sumoN &Y^*_{ik}\overline{V_k}[1] + V_i[1]\left(\sumoN Y^*_{ik}\overline{V_k}[0]-a^*_i\right)= \\
&S_i-V_i[0]a^*_i-b^*_i,\qquad\qquad \qquad\quad i \in \mc{B}_{PQ}
 \end{aligned}\label{eq:gen-powseries1}\end{equation}
 and for $n\geq 2$
\begin{equation}\begin{aligned}
\overline{V_i}[0]&\sumoN Y_{ik}V_k[n] + \overline{V_i}[k](\sumoN Y_{ik}V_k[0]-a_i)=\\
 &-\!\!\sum_{m=1}^{n-1}\!\overline{V_i}[m]\!\!\left(\sumoN\! Y_{ik}V_k[n-m]+ \delta_{m,n-1}a_i\!\right), i \in \mc{B}_{PQ}\\
V_i[0]&\sumoN Y^*_{ik}\overline{V_k}[n] + V_i[k](\sumoN Y^*_{ik}\overline{V_k}[0]-a^*_i)=\\
&-\!\!\sum_{m=1}^{n-1}\!V_i[m]\!\!\left(\sumoN\! Y_{ik}^*\overline{V_k}[n-m]+\delta_{m,n-1}a_i^*\!\right), i \in \mc{B}_{PQ}
 \end{aligned}\label{eq:gen-powseries2}\end{equation}
 \\
 Similarly, equating coefficients of $z^n$ in \eqref{eq:genPV1} yields for $n=1$
 \begin{equation}\begin{aligned}
 \overline{V_i}&[0]\sumoN Y_{ik}V_k[1] + \overline{V_i}[1]\left(\sumoN Y_{ik}V_k[0]-a_i\right) +{}\\
 &V_i[0]\sumoN Y_{ik}^* \overline{V_k}[1]+V_i[1]\left(\sumoN Y^*_{ik}\overline{V_k}[0]-a^*_i\right) = \\
&\quad 2P_i-\overline{V_i}[0]a_i-V_i[0]a_i^*-b_i-b^*_i,\qquad i\in \mc{B}_{PV}
\end{aligned}\end{equation}
and for $n\geq2$
\begin{equation}\begin{aligned}
\overline{V_i}&[0]\sumoN Y_{ik}V_k[n] + \overline{V_i}[n]\left(\sumoN Y_{ik}V_k[0]-a_i\right) + \\
&V_i[0]\sumoN Y_{ik}^* \overline{V_k}[n]+V_i[n]\left(\sumoN Y^*_{ik}\overline{V_k}[0]-a^*_i\right) = \\
 &-\!\!\sum_{m=1}^{n-1}\overline{V_i}[m]\left(\sumoN Y_{ik}V_k[n\!-\!m]+\!\delta_{m,n-1}a_i\!\right)\\
 &-\!\!\sum_{m=1}^{n-1}V_i[m]\left(\sumoN Y_{ik}^*\overline{V_k}[n\!-\!m]+\!\delta_{m,n-1}a_i^*\!\right)
 \end{aligned}\end{equation}
 \\
 Finally, equating coefficients of $z^n$ in \eqref{eq:genPV2} gives
 \begin{equation}\begin{aligned}
 \overline{V_i}[0]V_i[n]+\overline{V_i}[n]V_i[0]=&-\sum_{m=1}^{n-1}\overline{V_i}[m]V_i[m-n] + L_i^2[n]\\
 \end{aligned}\label{eq:genpowseriesend}\end{equation}

Equations \eqref{eq:gen-powseries1}--\eqref{eq:genpowseriesend} can be written in the following simple form, where the matrix on the left is independent of $n$.
\begin{equation}
\begin{bmatrix} A_{PQ_1}&A_{PQ_2}\\ A_{PV_1}&A_{PV_2}\\ A_{PQ_3}&A_{PQ_4}\\A_{PV_3}&A_{PV_4} \end{bmatrix}\begin{bmatrix}V[n]\\ \overline{V}[n]\end{bmatrix}=\begin{bmatrix}r_{PQ_1,n\!-\!1}\\r_{PV_1,n\!-\!1}\\r_{PQ_2,n\!-\!1}\\r_{PV_2,n\!-\!1}\end{bmatrix}\label{eq:genblock}
\end{equation}

 where
 \begin{align*}
 A_{PQ_1,{ij}}&= \overline{V_i}[0]Y_{ij}\\
 A_{PV_1,{ij}}&= \overline{V_i}[0]Y_{ij}+\delta_{_{i,j}}\left(\sumoN Y_{ik}^*\overline{V_k}[0]-a_i^*\right)\\
 A_{PQ_2,{ij}}&= \delta_{_{i,j}}\left(\sumoN Y_{ik}V_k[0]-a_i\right)\\
 A_{PV_2,{ij}}&= V_i[0]Y^*_{ij}+\delta_{_{i,j}}\left(\sumoN Y_{ik}V_k[0]-a_i\right)\\
 A_{PQ_3,{ij}}&= \delta_{_{i,j}}\left(\sumoN Y_{ik}^*\overline{V_k}[0]-a_i^*\right)\\
 A_{PV_3,{ij}}&= \delta_{_{i,j}}\overline{V_i}[0]\\
 A_{PQ_4,{ij}}&= V_i[0]Y_{ij}^*\\
 A_{PV_4,{ij}}&= \delta_{_{i,j}}V_i[0]
 \end{align*}
and
 \begin{align*}
r_{PQ_1,n\!-\!1,i}&= \delta_{_{n,1}}\left(S_i^*-\overline{V_i}[0]a_i-b_i\right)-{}\\
& \quad\!\!\!\sum\limits_{m=1}^{n-1}\!\overline{V_i}[m]\!\!\left(\sumoN\limits \!\!Y_{ik}V_k[n\!-\!m]\!+\!\delta_{_{m,n-1}}a_i\!\!\right)\\
r_{PV_1,n\!-\!1,i}&= \delta_{_{n,1}}\left(2P_i-\overline{V_i}[0]a_i-V_i[0]a_i^*-b_i-b_i^*\right)-{}\\
&\quad\!\!\!\sum\limits_{m=1}^{n-1}\!\overline{V_i}[m]\!\!\left(\sumoN\limits \!\!Y_{ik}V_k[n\!-\!m]\!+\!\delta_{_{m,n-1}}a_i\!\!\right)-{}\\
&\quad\!\!\!\sum\limits_{m=1}^{n-1}\!V_i[m]\!\!\left(\sumoN\limits \!\!Y^*_{ik}\overline{V_k}[n\!-\!m]\!+\!\delta_{_{m,n-1}}a_i^*\!\!\right)\\
r_{PQ_2,n\!-\!1,i}&= \delta_{_{n,1}}\left(S_i-V_i[0]a^*_i-b^*_i\right)-{}\\
&\quad\!\!\!\sum\limits_{m=1}^{n-1}\!V_i[m]\!\!\left(\sumoN\limits \!\!Y^*_{ik}\overline{V_k}[n\!-\!m]\!+\!\delta_{_{m,n-1}}a_i^*\!\!\right)\\
r_{PV_2,n\!-\!1,i}&= -\!\!\sum\limits_{m=1}^{n-1}\overline{V_i}[m]V_i[n-m]+L_i^2[n]
\end{align*}

Values for $V_i[0]$ and $\overline{V}_i[0]$ can be obtained by solving
(\ref{eq:genPQ})--(\ref{eq:genPV2}) at $z=0$ or equivalently equating
the constant terms (coefficient for $z^0$). Unfortunately this leads
to a nonlinear system of equations for the $V_i[0],\overline{V}_i[0]$
and a solution for general $a_i, b_i, L_i$ can not be derived easily.
However, for each of the specific models that we consider it is
possible to state a simple choice for $V_i[0],\overline{V}_i[0]$ that
furthermore satisfies $\overline{V}_i[0] = V^*_i[0]$. This choice will be
used as the seed, allowing us to calculate all further coefficients of
the power series from the recurrence (\ref{eq:genblock}).

\subsection{Holomorphicity and Reflecting Condition}
We now assume that we have a seed solution $(v, \overline{v})$, with $v_i = V_i[0], \overline{v}_i= \overline{V}_i[0]$, to (\ref{eq:genPQ})--(\ref{eq:genPV2}) at $z=0$ and that furthermore this seed satisfies $\overline{v} = v^*$.

Proving that this seed can be continued into holomorphic functions $V(z)$ and $\overline{V}(z)$ of $z$ that solve the general model follows a similar approach as in Section \ref{sec:HELM-holo}. Converting equations (\ref{eq:genPQ})--(\ref{eq:genPV2}) into functions to replace $f$ in \eqref{eq:ift}, we need to prove that $J$, defined in \eqref{eq:CIFTJ}, is non-singular at $(0,v,\overline{v})$. Examination of $J$ and $A$ from \eqref{eq:genblock} at $(0,v,\overline{v})$ shows that the two matrices are equivalent at this point. 
Therefore $A$ being
nonsingular is sufficient for $V(z),\overline{V}(z)$ to be holomorphic using the CIFT.
While it is possible to create special networks where this does not
hold, especially if no restrictions are placed on the values of $a$
and $b$, for general networks and sensibly chosen values for $a$ and
$b$ the condition will hold, as indeed is true for the IEEE models
used for the numerical tests in Section~\ref{sec:results}.

The proof that the reflecting condition \eqref{eq:helm-refcond} is
implied follows the same pattern as the proof given in Section \ref{sec:refcond1}.
Note that for the PQ-part of (\ref{eq:genblock}) we have
$r_{PQ_1,n\!-\!1}=r_{PQ_2,n\!-\!1}^*$. In the PV part,
$r_{PV_1,n\!-\!1}$ and $r_{PV_2,n\!-\!1}$ are both real-valued:
Assuming that we have already proven 
that $\overline{V}_i[k] = V^*_i[k]$ for all $k<n$, then $r_{PV_1,n\!-\!1}$, can be written as the sum of a complex number and its complex conjugate, whereas for $r_{PV_2,n\!-\!1}$ we have $\sum_{m=1}^{n-1}\overline{V_i}[m]V_i[n-m] \in \mc{R}$ since it is
equal to its own complex conjugate, and $L_i(z)\in \mc{R}$ by assumption. 

With this in mind, taking the complex conjugates of both sides of (\ref{eq:genblock}) we obtain
\begin{equation}
\begin{bmatrix} A_{PQ_1}^\ast&A_{PQ_2}^*\\ A_{PV_1}^*&A_{PV_2}^*\\A_{PQ_3}^*&A_{PQ_4}^*\\A_{PV_3}^*&A_{PV_4}^* \end{bmatrix}\begin{bmatrix}V[n]^*\\ \overline{V}[n]^*\end{bmatrix}=\begin{bmatrix}r_{PQ_1,n\!-\!1}^*\\r_{PV_1,n\!-\!1}^*\\r_{PQ_1,n\!-\!1}^*\\r_{PV_2,n\!-\!1}^*\end{bmatrix}
\end{equation}
which, after noting how $A_{PQ_1}^* = A_{PQ_4}$, $A_{PQ_2}^* = A_{PQ_3}$, $A_{PV_1}=A_{PV_2}$, and $A_{PV_3}=A_{PV_4}$, we can rearrange and rewrite as
\begin{equation}
\begin{bmatrix} A_{PQ_1}&A_{PQ_2}\\ A_{PV_1}&A_{PV_2}\\A_{PQ_3}&A_{PQ_4}\\A_{PV_3}&A_{PV_4} \end{bmatrix}\begin{bmatrix}\overline{V}[n]^*\\ V[n]^*\end{bmatrix}=\begin{bmatrix}r_{PQ_1,n\!-\!1}\\r_{PV_1,n\!-\!1}\\r_{PQ_2,n\!-\!1}\\r_{PV_2,n\!-\!1}\end{bmatrix}
\label{eq:refcondpv2}\end{equation}

Comparing \eqref{eq:refcondpv2} with the original system \eqref{eq:genblock} and assuming as above that $A$ is non-singular we must have that $V[n]=\overline{V}^*[n]$ and $V^*[n]=\overline{V}[n]$ for $n\geq 0$ as long as the reflecting condition holds for the seed $V_i[0] = V_i(0), \overline{V}_i[0] = \overline{V}_i(0)$. The reflecting condition $\overline{V}(z)=V(z^*)^*$ then follows.

Different choices of $a_i$, $b_i$, and $L_i(z)$ result in different paths for the bus voltages between $z\!=\!0$ and $z\!=\!1$. The obvious choices for $a_i$ and $b_i$ are either to both be 0, eliminating the terms from the model, or for one to be 0 and the other to be $y_i$, which makes the system of equations trivial at $z=0$. There is no single obvious choice for $L_i(z)$, though we have decided to have the voltage magnitudes scale linearly with $z$ in our models. 
Table \ref{tab:bustypes1} provides the values for $a_i, b_i,$ and $L_i(z)$ that we have chosen for the four models we present in this paper.
\begin{table}[!htbp]
\centering\footnotesize
\caption{Outline of PV methods}
\begin{tabular}{ccccc}
\toprule
Model& $a_i$&$b_i$&$L_i(z)$&$V_i(0)$\\
\midrule
1 & $y_i$ & 0 &$ 1+z(M_i-1) $& 1\\
2 & $\sumoN Y_{ik}\lambda_k$ & 0 & $M_i$ & $\lambda_i$\\
3 & 0 & 0 & $\|\nu_i\|+z(M_i-\|\nu_i\|)$ & $\nu_i$\\
4 & 0 & $y_i $& $1+z(M_i-1) $& 1\\
\bottomrule
\end{tabular}
\label{tab:bustypes1}
\end{table}

\subsection{Model 1}
\label{sec:mod4}
The first of our models sets $a_i=y_i$, $b_i=0$, and $L_i=1+z(M_i-1)$. With these values, we obtain a model that is very similar to the Subramanian model given in \cite{subramanian}, differing only in that the double convolutions are removed. Using the reflecting condition, the holomorphic embedding simplifies to
\begin{equation}\begin{aligned}
 \textrm{Re} &\left\{\!\!V_i^*(z^*)\left(\!\sumoN Y_{ik}V_k(z)-(1\!\!-\!\!z)y_i\right)\!\! \right\}=zP_i,&&i \in \mc{B}\\
 \textrm{Im} &\left\{\!\!V_i^*(z^*)\left(\!\sumoN Y_{ik}V_k(z)-(1\!\!-\!\!z)y_i\right)\!\! \right\}=-zQ_i,\!\!\!\!\!&& i \in \mc{B}_{PQ}\\
V_i^*&(z^*)V_i(z) = (1+z(M_i-1))^2,&& i\in\mc{B}_{PV}
 \end{aligned}\label{eq:mod4-hol}\end{equation}
The last equation gives the PV buses a voltage magnitude of 1 when at $z\!=\!0$, which allows the initial model to have the simple solution of $V_i(0)=1,i\in \mc{B}$ when $z\!=\!0$, while allowing the voltages in the PV buses to reach their required magnitude, $M_i$, when $z\!=\!1$. All the equations are now real, so we split the voltage $V_i$ into its real and imaginary parts, and so equation \eqref{eq:genblock} becomes 
\begin{equation}
\begin{bmatrix} A_{_1}&A_{_2}\\ A_{PQ_3}&A_{PQ_4}\\A_{PV_3}&A_{PV_4} \end{bmatrix}\begin{bmatrix}\textrm{Re}\{V[n]\}\\ \textrm{Im}\{V[n]\}\end{bmatrix}=\begin{bmatrix}\textrm{r}_{_{1,n-1}}\\\textrm{r}_{PQ_{2,n-1}}\\\textrm{r}_{PV_{2,n-1}}\end{bmatrix}\label{eq:mod1block}
\end{equation}
 
where the A matrix is
 \begin{align*}
A_{1_{ij}}&= G_{ij}, && i\in \mc{B}\\
 A_{2_{ij}}&= B_{ij}, && i\in \mc{B}\\
 A_{PQ_{3_{ij}}}&= B_{ij}, && i\in \mc{B}_{PQ}\\
 A_{PV_{3_{ij}}}& = 2\delta_{i,j},&& i\in\mc{B}_{PV}\\
 A_{PQ_{4_{ij}}}&=  G_{ij},&& i\in \mc{B}_{PQ}\\
  A_{PV_{4_{ij}}}&= 0,&& i\in\mc{B}_{PV}
 \end{align*}
 
and the right-hand side becomes

 \begin{align*}
 \!\!\textrm{r}&_{_{1,n-1,i}}= \delta_{_{n,1}}\left(P_i-\text{Re}\left\{y_i\right\}\right)-&{}\\
&\quad\!\!\!\text{Re}\left\{\!\sum\limits_{m=1}^{n-1}\!V_i^*[m]\!\!\left(\sumoN\limits \!\!Y_{ik}V_k[n\!-\!m]\!+\!\delta_{_{m,n-1}}y_i\!\!\right)\!\!\right\},&&i\in\mc{B}\\
 \!\!\textrm{r}&_{PQ_{2,n-1,i}}= \delta_{_{n,1}}\left(-Q_i-\text{Im}\left\{y_i\right\}\right)-&{}\\
&\quad\!\!\!\text{Im}\left\{\!\sum\limits_{m=1}^{n-1}\!V_i^*[m]\!\!\left(\sumoN\limits \!\!Y_{ik}V_k[n\!-\!m]\!+\!\delta_{_{m,n-1}}y_i\!\!\right)\!\!\right\},&&i\in\mc{B}_{PQ}\\
 \textrm{r}&_{PV_{2,n-1,i}}= -\!\!\sum\limits_{m=1}^{n-1}V_i^*[m]V_i[n-m]+L_i^2[n], &&i\in\mc{B}_{PV}
\end{align*}

We may calculate the power series of $V$ (and $\overline{V}$) to any desired degree $n$. This model does not suffer from the double convolution found in~\cite{subramanian}, and in Section \ref{sec:results} we will show that this model leads to more accurate solutions of the IEEE test cases.

\subsection{Model 2}
\label{sec:mod2}
The idea behind Model 2 is to keep the voltage magnitudes of PV buses constant regardless of the value of $z$. This is accomplished by first setting $a_i=\sumoN Y_{ik}\lambda_k$ and $b_i=0$ in \eqref{eq:genPQ} and \eqref{eq:genPV1}, where all $\lambda_k$ are constant and $\lambda_k=1, k\in\mc{B}_{PQ}\cup\{0\}$, and $\lambda_k=\|V_k\|, k\in\mc{B}_{PV}$. As well, in \eqref{eq:genPV2} we set $L_i(z)=\|V_i\|$, which is independent of $z$. The holomorphic embedding thus becomes:

\begin{equation}\begin{aligned}
\textrm{Re}&\left\{ V_i^*(z^*)\left(\sumoN Y_{ik}V_k(z)-{}\right.\right.\\
&\left.\left.(1\!\!-\!\!z)\sumoN Y_{ik}\lambda_k\right)\!\!\right\}=zP_i,\quad i \in \mc{B}\\
\textrm{Im}&\left\{ V_i^*(z^*)\left(\sumoN Y_{ik}V_k(z)-{}\right.\right.\\
&\left.\left.(1\!\!-\!\!z)\sumoN Y_{ik}\lambda_k\right)\!\!\right\}=-zQ_i,\quad i \in \mc{B}_{PQ}\\
 V_i^*&(z^*)V_i(z) = \lambda_i^2 ,\quad i\in\mc{B}_{PV}
\end{aligned}\label{eq:mod2-holo}\end{equation} 

For Model 2, when $z\!=\!0$, $V_i(0)=\lambda_i, \overline{V_i}(0)=\lambda_i$ is a valid solution for $i\in\mc{B}$. 

Here, the $A$ matrix from equation \eqref{eq:mod1block} becomes:
 \begin{align*}
A_{1_{ij}}&= G_{ij}\lambda_i, && i\in \mc{B}\\
A_{2_{ij}}&= B_{ij}\lambda_i, && i\in \mc{B}\\
A_{PQ_{3_{ij}}}&= B_{ij}\lambda_i, && i\in \mc{B}_{PQ}\\
A_{PV_{3_{ij}}}& = 2\delta_{i,j}\lambda_i,&& i\in\mc{B}_{PV}\\
A_{PQ_{4_{ij}}}&=  G_{ij}\lambda_i,&& i\in \mc{B}_{PQ}\\
A_{PV_{4_{ij}}}&= 0,&& i\in\mc{B}_{PV}
 \end{align*}
 
while the right-hand side becomes

 \begin{align*}
 \!\!\textrm{r}&_{_{1,n-1,i}}= \delta_{_{n,1}}\left(P_i-\lambda_i\sumoN G_{ik}\lambda_k\right)-&{}\\
&\ \!\!\text{Re}\!\left\{\!\sum\limits_{m=1}^{n-1}\!\!V_i^*[m]\!\!\left(\sumoN\limits \!\!Y_{ik}V_k[n\!-\!m]\!+\!\delta_{_{m,n\!-\!1}}\!\!\sumoN \!Y_{ik}\lambda_k\!\!\right)\!\!\!\right\}\!\!,\ \!i\in\mc{B}\\
 \!\!\textrm{r}&_{PQ_{2,n-1,i}}= \delta_{_{n,1}}\left(-Q_i-\lambda_i\sumoN B_{ik}\lambda_k\right)-&{}\\
&\ \!\!\text{Im}\!\left\{\!\sum\limits_{m=1}^{n-1}\!\!V_i^*[m]\!\!\left(\sumoN\limits \!\!Y_{ik}V_k[n\!-\!m]\!+\!\delta_{_{m,n\!-\!1}}\!\!\sumoN \!Y_{ik}\lambda_k\!\!\right)\!\!\!\right\}\!\!,\ \!i\in\mc{B}_{PQ}\\
 \textrm{r}&_{PV_{2,n-1,i}}= -\!\!\sum\limits_{m=1}^{n-1}V_i^*[m]V_i[n-m],\qquad i\in\mc{B}_{PV}
\end{align*}

The $A$ matrix is identical to that in Model 1 save the introduction of $\lambda$ -- the main difference is in the right-hand side. We can once again solve for $V$ to any desired value of $n$.
\subsection{Model 3}
\label{sec:mod3}
In Model 3 both $a_i$ and $b_i$ are set to 0. This gives a model similar in form to the model used in the HELM method given in \cite{helmpatent}---in the absence of PV buses they are identical. The solution at $z=0$ requires solving a simple set of equations. We obtain $V(0)=\nu$ and $\overline{V}(0)=\nu^*$, where $\nu$ is the solution to 

\begin{equation}
\sumoN Y_{ik} \nu_k = 0\qquad i\in\mc{B}
\end{equation}

We also scale the voltage magnitudes linearly with respect to $z$, so $L_i(z)= \|\nu_i\|+z(M_i-\|\nu_i\|)$. The holomorphic embedding thus becomes:
\begin{equation}\begin{aligned}
& \textrm{Re}\left\{V_i^*(z^*)\sumoN Y_{ik}V_k(z)\!\right\}=zP_i , && i \in \mc{B}\\
& \textrm{Im}\left\{V_i^*(z^*)\sumoN Y_{ik}V_k(z)\!\right\}=-zQ_i , && i \in \mc{B}_{PQ}\\
& V_i^*(z^*)V_i(z) = \left(\|\nu_i\|+z(M_i-\|\nu_i\|)\right)^2 , && i\in\mc{B}_{PV}
\end{aligned}\label{eq:mod3-holo}\end{equation} 

Here, the $A$ matrix from equation \eqref{eq:mod1block} becomes:
 \begin{align*}
A_{1_{ij}}&= \nu_iG_{ij}, && i\in \mc{B}\\
A_{2_{ij}}&= \nu_iB_{ij}, && i\in \mc{B}\\
A_{PQ_{3_{ij}}}&= \nu_iB_{ij}, && i\in \mc{B}_{PQ}\\
A_{PV_{3_{ij}}}& = 2\delta_{i,j}\nu_i,&& i\in\mc{B}_{PV}\\
A_{PQ_{4_{ij}}}&=  \nu_iG_{ij},&& i\in \mc{B}_{PQ}\\
A_{PV_{4_{ij}}}&= 0,&& i\in\mc{B}_{PV}
 \end{align*}
  
while the right-hand side becomes
 \begin{align*}
 &\textrm{r}_{_{1,n-1,i}}= \delta_{_{n,1}}P_i-\text{Re}\left\{\sum\limits_{m=1}^{n-1}\!V_i^*[m]\!\!\sumoN\limits \!\!Y_{ik}V_k[n\!-\!m]\!\right\},\quad\!\! i\in\mc{B}\\
 &\textrm{r}_{PQ_{2,n-1,i}}= -\delta_{_{n,1}}Q_i-{}\\
 &\qquad\qquad\qquad \text{Im}\left\{\sum\limits_{m=1}^{n-1}\!V_i^*[m]\!\!\sumoN\limits \!\!Y_{ik}V_k[n\!-\!m]\!\right\},\quad i\in\mc{B}_{PQ}\\
 &\textrm{r}_{PV_{2,n-1,i}}= -\!\!\sum\limits_{m=1}^{n-1}V_i^*[m]V_i[n-m]+L_i^2[n], \quad i\in\mc{B}_{PV}
\end{align*}

We can once again solve for $V$ to any desired value of $n$.
\subsection{Model 4}
The holomorphic embedding for Model 4 is the same as Model 1 except that the term involving $(1\!-\!z)$ is made voltage independent. This is accomplished by  setting $a_i=0$, $b_i=y_i$, and $L_i=1+z(M_i-1)$. The holomorphic embedding becomes
\begin{equation}\begin{aligned}
 & \textrm{Re}\left\{V_i^*(z^*)\sumoN \!\!Y_{ik}V_k(z)-(1\!\!-\!\!z)y_i\right\}=zP_i, &&\!\!\!i \in \mc{B}\\
 & \textrm{Im}\left\{V_i^*(z^*)\sumoN \!\!Y_{ik}V_k(z)-(1\!\!-\!\!z)y_i\right\}=-zQ_i,&&\!\!\!i \in \mc{B}_{PQ}\\
 &V_i^*(z^*)V_i(z) = (1+z(M_i-1))^2, &&\!\!\!i\in\mc{B}_{PV}
\end{aligned}\label{eq:mod5-holo}\end{equation} 

Here, the $A$ matrix from equation \eqref{eq:mod1block} becomes:
 \begin{align*}
A_{1_{ij}}&= G_{ij}+\delta_{_{i,j}}\textrm{Re}\left\{y_i\right\}, && i\in \mc{B}\\
A_{2_{ij}}&= B_{ij}-\delta_{_{i,j}}\textrm{Im}\left\{y_i\right\}, && i\in \mc{B}\\
A_{PQ_{3_{ij}}}&= B_{ij} -\delta_{_{i,j}}\textrm{Im}\left\{y_i\right\},&& i\in \mc{B}_{PQ}\\
A_{PV_{3_{ij}}}& = 2\delta_{i,j},&& i\in\mc{B}_{PV}\\
A_{PQ_{4_{ij}}}&=  G_{ij}+\delta_{_{i,j}}\textrm{Im}\left\{y_i\right\},&& i\in \mc{B}_{PQ}\\
A_{PV_{4_{ij}}}&= 0,&& i\in\mc{B}_{PV}
 \end{align*}
  
while the right-hand side becomes
 \begin{align*}
 \textrm{r}_{_{1,n-1,i}}&= \delta_{_{n,1}}\left(P_i-\text{Re}\left\{y_i\right\}\right)-&{}\\
&\quad\text{Re}\left\{\!\sum\limits_{m=1}^{n-1}\!V_i^*[m]\!\sumoN\limits \!\!Y_{ik}V_k[n\!-\!m]\!\!\right\},&&i\in\mc{B}\\
 \textrm{r}_{PQ_{2,n-1,i}}&= \delta_{_{n,1}}\left(-Q_i-\text{Im}\left\{y_i\right\}\right)-&{}\\
&\quad\text{Im}\left\{\!\sum\limits_{m=1}^{n-1}\!V_i^*[m]\!\sumoN\limits \!\!Y_{ik}V_k[n\!-\!m]\!\!\right\},&&i\in\mc{B}_{PQ}\\
 \textrm{r}_{PV_{2,n-1,i}}&= -\!\!\sum\limits_{m=1}^{n-1}V_i^*[m]V_i[n-m]+L_i^2[n], &&i\in\mc{B}_{PV}
\end{align*}

and we can once again solve for $V$ to any desired value of $n$.

Having shown how the new models are able to create the required power series for the bus voltages, the next section will show how successful each model was at solving IEEE test cases.

\section{Computational results}
\label{sec:results}
The models from Section~\ref{sec:PVmodels} were created in Matlab and run on a series of 
seven standard IEEE test cases (9-, 14-, 30-, 39-, 57-, 118-, 
and 300-Bus) to obtain power series for the voltages. 
These power series were then run through the Viskovatov Pad\'e approximant algorithm, as this is also the algorithm used in~\cite{helm}~and~\cite{subramanian}.
The resulting voltage values were then checked to see how well they solved the initial BPEE (\ref{eq:bpee}) as well as how closely they resemble the solution obtained through the power flow function in MatPower 4.1, used here as an example of a traditional method of solving a load flow problem.

In Table \ref{tab:res}, Max $|Rs|$ 
refers to the maximum absolute residual obtained from the 
$(N\!-\!1)$ equations when the Pad\'e approximant solutions are substituted back 
into (\ref{eq:bpee}). Max~$|\Delta|$ 
refers to the 
maximum absolute difference between the voltages obtained from the models and those obtained from MatPower. In each case the 
[15/15] Pad\'e approximant is used to determine the bus voltage values.  In all cases double precision was used.
The times given in the table are the average of 100 runs. All the models, including MatPower, take roughly the same amount of time for the small systems.

\begin{table}[htbp]
\centering\footnotesize
\caption{Results for IEEE test cases}
\begin{tabular}{rcccr}
\toprule
System & Model & Max |$R_s$| & Max $|\Delta|$ & Time (s)\\
\midrule
\multirow{6}{*}{case9} & Subr. 1 & 2.1932e-02 & 4.7824e-03 & 0.0110 \\
 & Subr. 2 & 1.2074e-11 & 1.4392e-12 & 0.0064 \\
 & 1 & 1.8938e-11 & 1.8491e-12 & 0.0059 \\
 & 2 & 1.8948e-11 & 1.8494e-12 & 0.0059 \\
 & 3 & 8.5691e-11 & 1.1914e-11 & 0.0073 \\
 & 4 & 4.4744e-12 & 6.1133e-13 & 0.0049 \\
 & MatPower & 5.9274e-14 && 0.0171 \\
\midrule
\multirow{6}{*}{case14} & Subr. 1 & 7.4505e-03 & 1.3738e-03 & 0.0111 \\
 & Subr. 2 & 2.3921e-14 & 5.8243e-12 & 0.0052 \\
 & 1 & 6.9950e-14 & 5.8213e-12 & 0.0062 \\
 & 2 & 1.1199e-12 & 5.6706e-12 & 0.0053 \\
 & 3 & 2.4484e-14 & 5.8254e-12 & 0.0049 \\
 & 4 & 2.4461e-14 & 5.8235e-12 & 0.0050 \\
 & MatPower & 1.3180e-10 && 0.0087 \\
\midrule
\multirow{6}{*}{case30} & Subr. 1 & 1.6959e-02 & 1.5330e-03 & 0.0147 \\
 & Subr. 2 & 2.2901e-14 & 1.9658e-10 & 0.0062 \\
 & 1 & 4.8474e-14 & 1.9658e-10 & 0.0068 \\
 & 2 & 2.3742e-14 & 1.9658e-10 & 0.0062 \\
 & 3 & 4.6547e-14 & 1.9658e-10 & 0.0059 \\
 & 4 & 6.0382e-14 & 1.9658e-10 & 0.0060 \\
 & MatPower & 9.7323e-10 && 0.0091 \\
\midrule
\multirow{6}{*}{case39} & Subr. 1 & 1.2063e+01 & 8.4773e-01 & 0.0198 \\
 & Subr. 2 & 2.2529e-01 & 2.5868e-02 & 0.0076 \\
 & 1 & 7.5746e-07 & 4.0769e-08 & 0.0085 \\
 & 2 & 5.0394e-06 & 2.4831e-07 & 0.0075 \\
 & 3 & 3.0385e-06 & 9.7468e-08 & 0.0069 \\
 & 4 & 1.1003e-09 & 5.2491e-11 & 0.0067 \\
 & MatPower & 8.2567e-13 && 0.0081 \\
\midrule
\multirow{6}{*}{case57} & Subr. 1 & 5.8562e+00 & 1.8105e-01 & 0.0181 \\
 & Subr. 2 & 4.7931e-13 & 2.2841e-13 & 0.0084 \\
 & 1 & 4.0425e-10 & 8.4430e-11 & 0.0092 \\
 & 2 & 2.7227e-09 & 1.0797e-10 & 0.0092 \\
 & 3 & 2.4653e-10 & 1.2281e-10 & 0.0088 \\
 & 4 & 4.8125e-10 & 2.7309e-10 & 0.0089 \\
 & MatPower & 3.7036e-12 && 0.0094 \\
\midrule
\multirow{6}{*}{case118} & Subr. 1 & 2.2228e+02 & 2.4416e+00 & 0.0713 \\
 & Subr. 2 & 8.1703e-11 & 4.3043e-12 & 0.0232 \\
 & 1 & 1.3121e-04 & 3.1500e-06 & 0.0192 \\
 & 2 & 2.0368e-04 & 1.6476e-05 & 0.0298 \\
 & 3 & 2.6099e-02 & 2.0917e-04 & 0.0220 \\
 & 4 & 1.6917e-10 & 7.6155e-12 & 0.0174 \\
 & MatPower & 1.4892e-12 && 0.0111 \\
\midrule
\multirow{6}{*}{case300} & Subr. 1 & 1.9832e+02 & 9.1111e+00 & 0.1348 \\
 & Subr. 2 & 1.4294e-01 & 4.9324e-02 & 0.0807 \\
 & 1 & 7.8944e+01 & 2.7515e+00 & 0.0656 \\
 & 2 & 1.5028e+02 & 2.4886e+00 & 0.1225 \\
 & 3 & 5.7531e+04 & 2.7865e+01 & 0.0761 \\
 & 4 & 2.8486e-04 & 8.4840e-06 & 0.0656 \\
 & MatPower & 1.7628e-12 && 0.0229 \\
\bottomrule
\end{tabular}
\label{tab:res}
\end{table}

From Table \ref{tab:res} we can see that the model derived from~\cite{subramanian}, Subr.~1, converges much more poorly than 
the rest. For the 9-Bus system, Subramanian's first model required a [60/60] Pad\'e 
approximant to achieve a similar level of accuracy to what the other models 
achieve with the [15/15] Pad\'e approximant. 
The other models perform relatively equally on all systems up to the 300-Bus system except that Subramanian's second model has difficulty with the 39-Bus system. All of the models have trouble with the 300-Bus system and fail to provide correct voltages except Model 4, which manages to converge slowly to a proper solution.

\begin{figure}[htbp]
\centering\footnotesize
\input{case300-Model5-VC.tikz}
\vspace{-1em}
\caption{Power Series Coefficients using Model 4 for 300-bus network.}
\label{fig:mod5vc}
\end{figure}
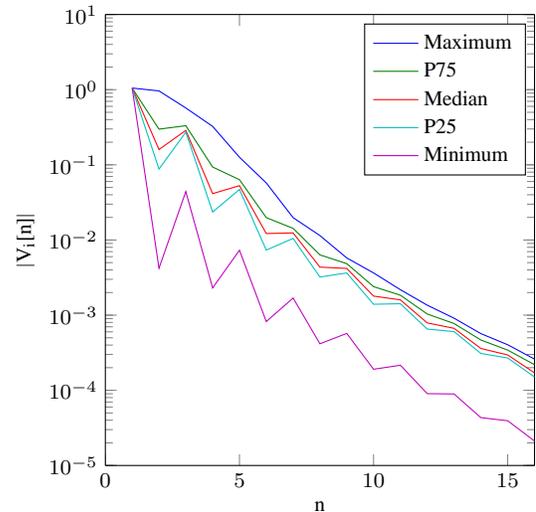

Below we will further investigate the behaviour for the 
300-bus system.
Figures~\ref{fig:mod4vc}~and~\ref{fig:mod5vc} show the power series 
coefficients obtained in the 300-Bus network using Models~1~and~4 respectively, while 
Figures~\ref{fig:mod4v1}~and~\ref{fig:mod4v2} show the corresponding singularities (zeros of the denominator polynomial) of the [50/50]
Pad\'e approximant. The data for the figures was produced using Maple to 100 significant digits.

We can see from Figure~\ref{fig:mod4v1}
that the Bus 8 Pad\'e approximant derived using Model~1 appears to have a set of poles along the real line between $z=0.1$ and $z=0.2$, and going off the real axis into the upper and lower half-planes, with yet another set around $z=0.8$, indicating cuts in the complex voltage function. This is consistent with a small radius of convergence of the corresponding power series as can be seen from its rapidly increasing coefficients (Figure~\ref{fig:mod4vc}). 
From Table \ref{tab:res} we see that these singularities adversely affect the rate of convergence of the Pad\'e approximants, and using double precision the model is incapable of providing a sufficiently accurate value at $z=1$. 
Further investigation was conducted on the 300-Bus network using Maple, where precision can be set to higher levels than double precision. By setting the precision to 200 digits, much higher order of Pad\'e approximants can be computed to greater accuracy, and in this setup convergence to the correct value at $z=1$ is obtained.

On the other hand, Figure~\ref{fig:mod4v2}, 
shows for Model~4 there are no singularities in the disk centred at $z=0$ of radius 1. The single pole closest to the origin below the real line is a spurious pole and does not appear neighbouring Pad\'e Approximants. Indeed the radius of convergence of the power series is greater than 1 as is indicated by the decreasing nature of its coefficients (Figure~\ref{fig:mod5vc}).

\begin{figure}[htbp]
\centering\footnotesize
\input{case300-Bus8-Model1Singularities.tikz}
\vspace{-1em}
\caption{Pad\'e Singularities for Bus 8 using Model 1}
\label{fig:mod4v1}
\end{figure}
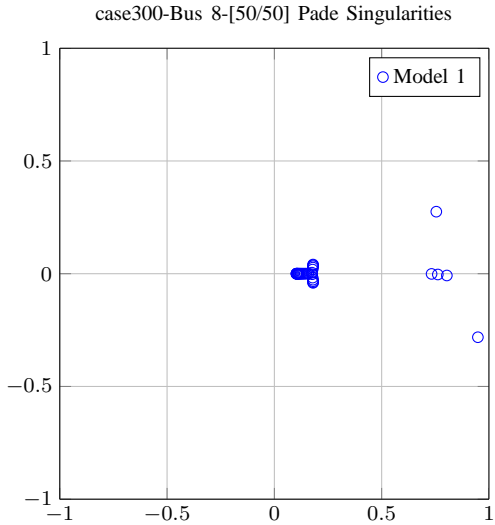
\begin{figure}[htbp]
\centering\footnotesize
\input{case300-Bus8-Model4Singularities.tikz}
\vspace{-1em}
\caption{Pad\'e Singularities for Bus 8 using Model 4}
\label{fig:mod4v2}
\end{figure}
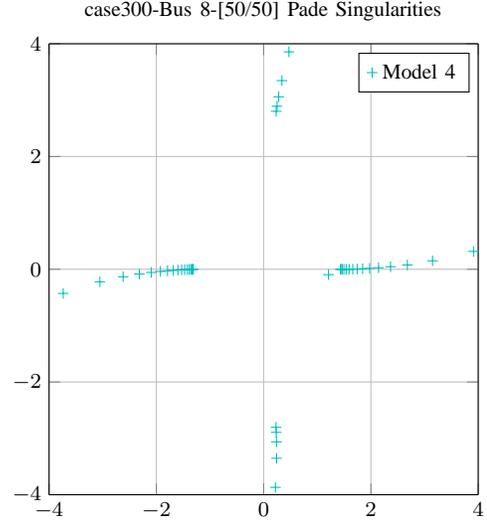

\begin{figure}[htbp]
\centering\footnotesize
\input{case300-Model4-VC.tikz}
\vspace{-1em}
\caption{Power Series Coefficients using Model 1 for 300-bus network.}
\label{fig:mod4vc}
\end{figure}
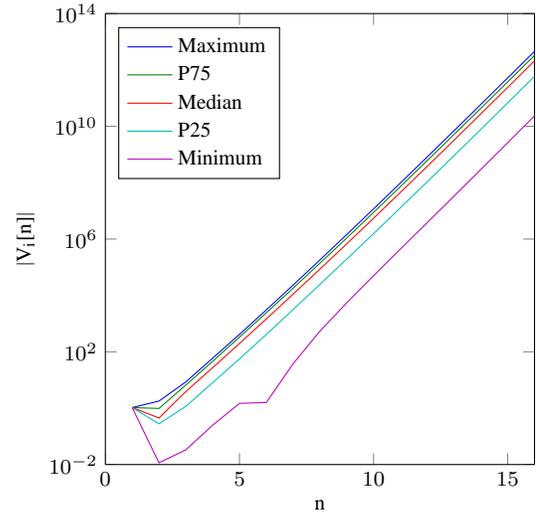

\subsection{Case 39PQ with HELM}
It is not possible to test the original HELM model from \cite{helm} on examples containing PV buses. Practical examples, including all of the IEEE test cases used by our general models, have PV buses.  One workaround to the lack of PQ-only networks is to convert the PV buses in the IEEE test cases into PQ buses. This is accomplished by finding the reactive power output of the PV buses at the power flow solution to the networks at full-load (as obtained by MatPower) and then fixing this value as a parameter at the bus while also now allowing voltage magnitude to vary (as in a PQ bus). The test-case solution should again be a solution to this altered network, and indeed, running these altered IEEE test cases through MatPower results in the same power flows and voltages as the original cases.

Using the HELM model on these modified results in the same
solution as MatPower in all cases except the 39-Bus IEEE
test case. In this case, the HELM model converges to a
solution at $z=1$ which does indeed solve the BPEE correctly, but which however
is different from the solution obtained by MatPower. The MatPower solution has voltages magnitudes close to 1 p.u. whereas the voltages magnitudes in the HELM solution are quite high (1.4--1.8 p.u.)

\begin{figure}[htbp]
\centering\footnotesize
\input{case39PQ-MatPowerBusVoltages.tikz}
\vspace{-1em}
\caption{Matpower Solutions going from no-load to full-load}
\label{fig:Matpower}
\end{figure}
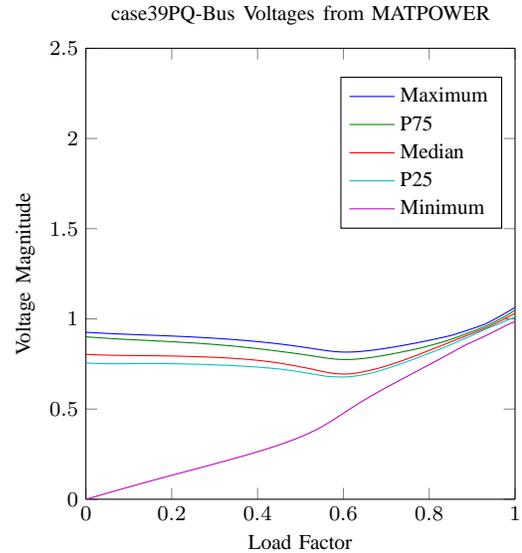

We have attempted to follow both solutions as the load changes between full load ($z=1$) and no load ($z=0$).  Figure~\ref{fig:39PQHelm}
shows the profile of the voltages at each bus as HELM moves from the initial
no-load case to the full-load case. We have attempted to use HELM to trace the MatPower solution back from $z=1$ to $z=0$ but have been unable to do so: there seems to be a singularity around $z=0.5$ which prevents HELM obtaining solutions for real values of $z$ less than $0.5$ due to precision problems.  Instead we have followed the MatPower solution back to the no-load case by using a traditional real homotopy on the non-embedded system (\ref{eq:bpeez}). The resulting bus voltage profiles as the load factor changes are shown in Figure~\ref{fig:Matpower}.

By comparing the two figures we see that the HELM solution starts with $V_i(0)\ne 0, \forall i$, whereas for the MatPower solution for one of the buses (Bus 38) $V_i(0)=0$.

The apparent inconsistencies can be reconciled by looking at the stability of each solution. Using the $\frac{d\Delta Q}{dV}$ stability criteria \cite{powersystems}, the MatPower solution---which is stable for the original IEEE test case---becomes unstable in the altered, PQ-only network. Indeed, if more reactive power is demanded at any PQ bus except Bus 32, then the voltage magnitude at that bus increases rather than decreases. The higher-voltage HELM solution, by contrast, is stable, giving credence to the claim that HELM will return a stable solution.

\begin{figure}[htbp]
\centering\footnotesize
\input{case39PQ-HELMBusVoltages.tikz}
\vspace{-1em}
\caption{HELM Solutions going from no-load to full-load}
\label{fig:39PQHelm}
\end{figure}
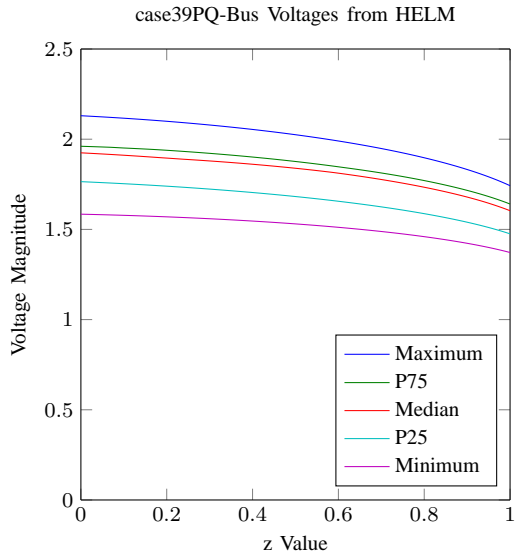

\bibliography{helm_pv_paper}

\end{document}

%% file: case300-Bus8-Model1Singularities.tikz
%
%
\definecolor{mycolor1}{rgb}{0.00000,0.75000,0.75000}%
\begin{tikzpicture}

\begin{axis}[%
width=5.705521cm,
height=6cm,
at={(0cm,0cm)},
scale only axis,
separate axis lines,
every outer x axis line/.append style={black},
every x tick label/.append style={font=\color{black}},
xmin=-1,
xmax=1,
xmajorgrids,
every outer y axis line/.append style={black},
every y tick label/.append style={font=\color{black}},
ymin=-1,
ymax=1,
ymajorgrids,
title={case300-Bus 8-[50/50] Pade Singularities},
legend style={legend cell align=left,align=left,draw=black}
]
\addplot [color=blue,only marks,mark=o,mark options={solid}]
  table[row sep=crcr]{%
.1026175257558457149180535546550495992171549750245793965437528786900505525156930196981611832204294875 +3.322243284374517701893684415044268586039430997900479112577119191039197952081563664319540392766503392e-7\\
.1029316869702344511143968536871391787386622352402129451069580415585436967314582816934409722456898766 +0.1342220988404814561594565945165161990124124405839185941861960780907701659924119874910381402733937553e-5\\
.1034583518639032764663965510847950179899733189803927731604220388078551434899279742309574643714848366 +0.3070611204117374776433230686907712141354035213167152436631072180870629122334537533075501037530489552e-5\\
.1042021926975519214624245103234610819169371948842783480474857342851137741074204962183787067374171839 +0.5587270516036863990428273304385308607129288173868008651790050047465333906472645014019065213766547564e-5\\
.1051699130718061263788889048578986842957894729351039858285361484011335051722930781693216105846511063 +0.8994606168347649501243319162767477064772889215646394107513556604472064107037594834483834259366822907e-5\\ 
.1063704269653290134572394913195839583864028593038996821277559324158596898682829842473931988527467978 +0.1343217021964435241463846072908096439976781294665077754786123051689648133041450534609341364598171442e-4\\ 
.1078151165363297002838562926985621502274743397596500258401251980178628790418931107849919069392742710 +0.1908264129827976991932610348092214303855290043644139505369731184485167778526094310485551489286000777e-4\\ 
.1095181971958619895398391990791547706953224586252160639834654145608262435218615814629136208427003715 +0.2617933084035213167273305347525464102470455956799251649441956106882951420431142991566095498020502051e-4\\ 
.1114972356915779271725076900006459965891054092832300126734349700085980756048625307417286816365932603 +0.3501581497772757602252978702788745300863108080551311979988523727485131212727701107356439643977224290e-4\\ 
.1137738948082484441627994291685672767495902338178491922516096734774555499153008566455677573418761584 +0.4595948148154436163054005031187075963402541222784019344160945758993769397747507003997719312670696389e-4\\ 
.1163750257502360249482580021640574591325846120284058577519723292987099365298683432623056481680219769 +0.5947338371244573737290680605258828036685269903589291972260728086886300008924494014111150803456566497e-4\\ 
.1193343173257695487591644154194554552920001435318750194341686233592867066043677851426748980104336459 +0.7615582088315533140087334928289575230406627293945400278577123543741682122974772361562896658830929080e-4\\ 
.1226948912165649868420617190805529654960108510050580195620308631520891793123031072899441404231411493 +0.9681656184456309753690146552784563916393457123800888850148128863929881163043591312420355732340730925e-4\\ 
.1265136333920683794869761290219881540309013698760819052898372807256950875801232154161468983887021245 +0.1226287996343785272071099514511849037512090509262319415660274606637343012125497411976278356629108324e-3\\ 
.1308690159706763291803009090629816393926978568440373641611764776231266356238948150380595852328448279 +0.1554479162239055794586544230188935029868509185085886581525149158564458413608344573599359401083468872e-3\\ 
.1358767413018797690609436883986949817252200843431189341267676195239748759030460186337353222707838473 +0.1985482465859728298374527969680102116220726399539240496285046496442226573062209967754091923740358079e-3\\ 
.1417256173490990437127252597456081545433505964198206696670211507498292983085689274559988825351356515 +0.2586232331149893409679421375456650454265498457123022559533119509450647501610647251321873979477464186e-3\\ 
.1487783996604590898654073107705363267668989335424775797966369901384731064607262039559410086166302633 +0.3528198321365088503105069458346223905992067497836844613431561898920675599319367718043534967181824326e-3\\ 
.1579818706209794991338676089291286209346811448073173436427062081719444659384438228471759259596781634 +0.5491524737653931240158090694991406342700958757726630576544274657298019975248517368767088402348431981e-3\\ 
.1750381238189407204448319746796062032079695384292581729031355392967121480587996568337217394473879774 +0.4641895484455543490618093410424060951819874953379964052172977236962344901547731256580649810511317104e-2\\ 
.1782131743657566491403673548525743815558517765009191898889511146677670865782558494194264617849775600 +0.1833589195179859569638658362248109750667283183972092616043940211713391367221977841456817781766146028e-1\\ 
.1798860433157709009217906820822707554656846600088249303845854720054502372306373664872422678420811765 +0.2709634378891019085303366903671204533580638718681784007838704034347085613169261934202267910543916735e-1\\ 
.1804578223843800192704027436411391917783634431911517107351985111349360321185642028197471967428143398 +0.3355822784502256788419058977387368230840292159937316336615174111883493268875587172702685230102126480e-1\\ 
.1804692396432930347982474151746447470351364575278465356944253229725130627667802848864706530004760567 +0.3800811728832671601570183320324961275973338482442550077600431291874721610488979923417715554475203286e-1\\ 
.1802806749195507106017216954372401686030905640913929360927789188214285631020900009988918631370815381 +0.4061178815969250050807553183551220564539284295491390671771413917314620133843590905038776498257423959e-1\\ 
.7546241143233838310449349810784500035640190122776773530473370586990742772704078538683232026917992338 +.2751026451751048004923897767265747480062487392818975536872530252318233277483098464818059490545145932\\ 
-3.302620668461887849770764268450595032621237888976305856183130754568262869924423029712554667255560353 +1.187055863643034011103249084787760389209069838972485058246484917952483793182229265737491061888341675\\ 
-1.718849939358709678264466229732453954587087464317662452908877927581001868013600693408604324717191404 +.3587415462565989394765942600611949995323105806067433112358903903220459755766482848635420319423992901\\ 
-1.873189477610871568801608576157476784892934587379027198211925082464232032473328947627619424298954607 +.3266882950113835000154245497363444712300320187310281896920335784492852172859644274796023727068819296\\ 
-1.016102898482174963308964145005438379130725511947421140828296845016841254374353116212808026654315081 +0.2453468135267641509128487184054955567289304024839995582156884736108115332259408760760554718492845816e-1\\ 
-1.631092181098703114641222531276241150735114100231288634071316067690612476665873227978093758224457183 -0.3249708532646971678922538577749100103761796003612941410324771809323605771344691559862459643825883198e-2\\ 
-1.890996782447336602564443721005135153476577324872159116189007969270821431690927174585757082178279556 -.3153854911150922199379885593179317598801055020687576068461766632185303525114289275179497944550139157\\ 
-7.808369987519317506055116462727990419013734927165314305535335859275994167706339635051007220224370008 -1.472870904092367811212546883866225186851757996855327440192190144867160224794300930876699476120682992\\ 
-1.721527612809177851886786285798734855807707701736021315371167131145797611521192873930134417579968781 -.3558974370756386242791643622561251157652560193811620475395375716555087982979471579334741982616850609\\ 
-3.165634100259231685988406846237335028966839198318469684441557645737873107355897223121288291046560165 -1.058727684747125574407754185408459717304998020132684772023742518737081053099186323468146616981872068\\ 
12.88979587786767610066179533522274220584089699823512080415493331955420181678141098725762424264850277 -14.39375790814018968923116892392593441485662450672096724287511808323626158526392028191990159044204592\\ 
6.797890631072776423367851922014275551358328613741913014383170269733833245279440689852499130417088899 -2.157677829485998343222541984857830905068037362225820953987945175476980917693973313380893805981672770\\ 
.9484187710387103319764039715363057624686252632346164244221723001047930978534059191172639662260618599 -.2817060100206696296133824596303933217102508888072835071883784355297323371724458340102174517196974566\\ 
.1802948663589302295964826954124776384414204546476778925905094003906717000282821046459810874581977308 -0.4057574823471415149237516987905598343593063606552439684820457372425498336248779971624093543076604192e-1\\ 
.1805305871068271582706950570658164124979563843929877849556014286314070512779092395332026362094497224 -0.3786610194309824695717473770855742670738923149680603115008030505015857507607341766653381555472035746e-1\\ 
.1805863474286578543495878444676939211047527108252252498049908916729913702597599997782221370996730531 -0.3324566907122255714244840257362123064142247793413816553793222912736444253890913393922577484584164206e-1\\ 
.1800705593504272454247335246488046587506795193875449649823873346727149987773605722071741333017681528 -0.2652900077944055472000246596084013160400119632374915067110544139576822208441541045192501149331054317e-1\\ 
.1784167760751797355160085862857901187592678083262359676406354196204585477214069897523239327847780038 -0.1736209816430897221161876749931989056915667605762180587579694717496448121208640380702968522135483386e-1\\ 
3.577599160362709516768862371249497881859037795862113015176739163896488101499613932590012046443643182 -.2063993225660228358787623474180867251688465492509712904084430720072279579531760180797662192597595834\\ 
1.508594502858655505306883662352436381141221903666791387717035760579119852966941224869712831724359440 -0.5458891736777488453899834129772341922374625685283238929347492146626580156700428641015179214050791449e-1\\ 
1.896596711152097504449610473740117024006181579377090064239658023051566971513748703934433275139357549 -0.6721537213742961493804807808449168637759729132818771079737149596410992007942571714795145525956631983e-1\\ 
.1762994363854823694832368007724382599841731686423330394976022976360431118043772981801706556108422282 -0.2381690765757164874404976567463088101686403111250751556931913993502943794454722074791544658230272301e-2\\ 
.8037545989419787002100420819704523427536918027312990701663147192490604176229169547045681921044257498 -0.7866029153773675177149673955683792795392604961603242254178163628044580431540214053383516709621406791e-2\\ 
.7617272736739866629682456075484740995867504174370269738822634263961503917003127628481095049571810878 -0.3444639851453956235419678484053028326118815321113847698410098454642757200683819618924984677195912343e-2\\ 
.7326077590574747532164462120113864363615455365612308036206419640526419551075223826261895110492368796 -0.8810793012214463810778548798742507550465808685237926037525673820548691841772203576340824624360110712e-3\\
};
\addlegendentry{Model 1};

\end{axis}
\end{tikzpicture}%

%% file: case300-Bus8-Model4Singularities.tikz
%
%
\definecolor{mycolor1}{rgb}{0.00000,0.75000,0.75000}%
\begin{tikzpicture}

\begin{axis}[%
width=5.705521cm,
height=6cm,
at={(0cm,0cm)},
scale only axis,
separate axis lines,
every outer x axis line/.append style={black},
every x tick label/.append style={font=\color{black}},
xmin=-4,
xmax=4,
xmajorgrids,
every outer y axis line/.append style={black},
every y tick label/.append style={font=\color{black}},
ymin=-4,
ymax=4,
ymajorgrids,
title={case300-Bus 8-[50/50] Pade Singularities},
legend style={legend cell align=left,align=left,draw=black}
]

\addplot [color=mycolor1,only marks,mark=+,mark options={solid}]
  table[row sep=crcr]{%
1.430218727913142819862804289093352227738965635707689267232716874125963154453218347078567326047155020 -0.6358676964813103183978662457912363936363172865596414422553912808153713948677082405928975141159193539e-4\\ 
1.661395171669803314524589357319781561464970357443732610443923683651731840278645814950832065650607759 +0.1178872865290438288227819170806188812650495607393431901045862336238796177051521754462425300296192340e-2\\ 
1.743099486610650092404212824665221767624942100143631012873751272841173860081181457281551905373621108 +0.4255735132170080974269972130719352110300453595659680673818058119864061623544820246106295681684936527e-2\\ 
1.845400178079268537995138514028307997776381029876627573737351333337685168807530552584687873276348046 +0.8751428638680835319379434635023387028238209568479931582491136617795077973935538852472766213862300698e-2\\ 
1.975658852219024575185223127951414736142035916075653989768208974468245464296190513662370743252621155 +0.1585227027933840045248338097392978673949869242502270072157017714809176510817083292963095327718068111e-1\\ 
2.142781122611890548099031960032927470491480302177155495065991962849993720142970699015615939745329831 +0.2709055394870155258381044923486096374478125939188966638262105727899459327563327573811946797542301609e-1\\ 
2.366117382475958355928131367977139084043210947101769968130346252624495067648885717006140213692056272 +0.4457933927227744984143000887013513635571172922735130933811338949373647138100830341408461861391515625e-1\\ 
2.678908034908888194731131395273424605616530746560783089281122382040300368565361957300226003094501149 +0.7683738066036304127251362472905035788969231051598277003604130969109564794560977914244447727029377548e-1\\ 
3.149023391789393809924351709652038597139334463421195553791341846177846236605200259182714202213327571 +.1497134609445888998843782046277241488857605821926259961724905408556912149984712374370785329177336325\\ 
3.913119491916982505850217304429973206780837223360215281916037888092253389229448326636386252207336685 +.3161479067945311465796795383902582861354124118147732714390430040995489681457530231298669169090466221\\ 
5.825085565220111196810483679416211095164611812923802047364266457405454818178565764784030396444610441 +1.018873219335658554785628812435214945847801461128238887255644950035982072455520917517490371134909437\\ 
4.094264046432934114560512238210038917350903875298286420706152887003969007338209946796434314057170539 +7.028407308485196826675920819296521617974033159978140342645003358072351924737995181419434192579773133\\ 
.8878571404845819010528013557609215066273609638220150771579002127581075871409535644020008281155736742 +4.925265088702839170679735204760785013307474841665269807029404034395554867597599588791396486664069189\\ 
.4712303944612762044146838613021364000991497424940395470588667369796103699673628744133155310699249510 +3.856015734211487946866750552524729496172765307331638835950872516828460487325352043740789274396982677\\ 
.3401939837424157407770342061041643039237823063422713771536346593570750986660206628432668064948167275 +3.346159074786780965115782740607220396167273009511924913093964597879809981789835615410882995511638872\\ 
.2803693939237443716240959809537257971632034311658176936678932306685328956914238112912239942188235862 +3.059298912884321949495267853681605214222420026719958588644230668382622551016424978684729806188016491\\ 
.2495064792909424638770709219912911758514610857168455358420541055488311627263556700417653193694197038 +2.891952292067629106111971443336918216230874146975450280564190183295174417846898412215159702295180533\\ 
.2342197905524246513169335268037513947679847837383287206755660973905090851928808358210686209411567903 +2.802953842490749136979092966081456869787097786043197743452008474773285962687108749470281344074546932\\ 
-1.310536609378397396379555889197881978419533075391224393813890869358429427312948185196743989739299941 -0.9805510998673405890428619820844767731101460503722011496225602421136921703509173462672208641086831457e-4\\ 
-1.319534452761198339334456827859588791456687866044167370456370430750085767760855055349360413791593210 -0.4004528494430370331044489283526878136941526398255623965424622288663814814286733172960613161033884656e-3\\ 
-1.334782003446976578797295723195966633291633214993262098053579901226168502479432220165400188060349183 -0.9313849052570937086898973658423559183868648860110142237627313232733950820647031372694235048991495992e-3\\ 
-1.356682673756577312886247298202259435334477749249244760462076114546163826794470477695667545425679489 -0.1732556691835775804708249295461830928179520062695111391916467818919025468963124905393688104593378818e-2\\ 
-1.385844149659852562472619032365782287548266334723783259736459385437166513895675042249454692067965182 -0.2868815804756849584512222088469190945220187056848087141907312127603561420333251708082858651073432691e-2\\ 
-1.423124215830241618705994297722682319721698430719641522441741568161205285038128082694489548015022339 -0.4438174975644591990411659040893875869913092855913325170665907507071440631940372752207098137345903146e-2\\ 
-1.469702754147566754121057372983445215740396860944244149420810785839641703906106819038276224727203751 -0.6589570958382172075052259848139205312049868483847112435119015313230340797389116169112396551045975928e-2\\ 
-1.527192166252638181509701363767553595174619102822636594748082900397325853415293535974737363774272067 -0.9564448997424396145492466020307835182384039913514659341418596145009692130505417373730722327109539953e-2\\ 
-1.597592968910641585699770325301317250658926037173578117341921512451189399974252464027683357342861642 -0.1370131886655049557163210764455349385606924316051275420798301716227646755734761269763552672822928399e-1\\ 
-1.684148902869247131897833819466568120121630907200936751708773628473400735083773747358362180799549467 -0.1934545610617884828776610187004528002727971519016630218602480705977195445971627753150056505732992798e-1\\ 
-1.791197266251156499945953797349744596861612184430777281685144881120563361756188768562844010176718752 -0.2748314350430931429890573218422137920648233070523740317911610706961921782653689030032784087858244884e-1\\ 
-1.924717692285124047843329411352553523811971583825455666936540220554713792983509470496220832268526351 -0.3922167327342752422664859760706569607256940992136814331680121643786414957853513008906593069472383644e-1\\ 
-2.094973515698846680684176325078643685325629852352156875792616262917816818145437880221712038267843450 -0.5683283655877572460690593868068782541828454830530589856795765653047809777830811420242237923530693102e-1\\ 
-2.317446093214470560264200529308299133791885221487742694412748821464219473481521530579628183172158843 -0.8504155179727333978886022311189188457408196474381574678415981714962817332264381799314545567109112048e-1\\ 
-2.619018176071538738157044095695407085663930536771918357160145537940343242720470114824769282089683935 -.1328700903840728828677631893081339578564708782100094943314811149854806004702164277210688472740750674\\ 
-3.053212109093833727315332102062405291471679674579053845622710240097241146074168027614525607346295843 -.2229400381734804698226271391254056729174817597073699264654808166368872364697300287658883819013286866\\ 
-3.736179428655587686524326130109731393078586143376752356357392177198357990730955684327476893494022899 -.4294617081144498985672753828528886437171066616081114605550544016384932787622888333515887933666170132\\ 
-4.936252606151495229397709069233093095043915664120232616035136870987608372781062274078999169608487398 -1.080442767902059918962815400173578580131858076085461705559301736251313354099956276211765975831995015\\ 
-6.438549961874290908238553842990585204704603133056845566792270057838050934316080648759005344068619989 -4.038901938203565159152431225039539907325615726134606871384325726324123430527431065927475295021745803\\ 
-1.952711792713595039452895065593922124116142363724459629826660748951674918946506470224269794669450429 -7.187774211377847630263147032758034409653372921640834117091347049752943121071388540635298968574365254\\ 
0.2805775786868690471320334127057226084521507740329343596065022971929712848052599478136580642414576338e-1 -4.926367766330307554253866010275649987287312716947136752657436832124807585528797658038563541746714075\\ 
.2201684029587129796918650886978670550016432317557569469471733122737529991505358498941801503822775390 -3.869157575245115500916132627386257577195856865554964416947854387971039095734995766278480220590460999\\ 
.2421390349874311181136737557455519764222466664113440143255819202640648150362679357593257129989949197 -3.352978130907295907141515132188181766761835873429217429335332359477945190123353736938559012469334373\\ 
.2398624337608398798957495813581644446139288842786207828838302595561349572011138281473315889243799902 -3.062391558020530973209912399694203867724593611770575836336515091369828752785076184426857370090348245\\ 
.2346520999665486103844674561412940350680298787857160529533113354432473493015886252101157940065659297 -2.893132534523735038881582345687977258398106049131493784662464064710971573532138163281119170097749354\\ 
.2308807110036017009410844092503745811454084806186535577150998373128511806773369612338950983503614438 -2.803223904590272442080701314150072896145067176528642956858865633929690107959706849010314591564840399\\ 
1.208789841802448164699006675744228418973667668821367392677221011764957572505098598894249380026355555 -0.9544148594329294792976551100988598438904711241024043019750559095931052489226206874589301578711028539e-1\\ 
1.539177274267164419772074244192025083088413308532743866552003767003366990976906112344077581651251965 -0.1003011566270993517078372285944734317226762552188120588705870827806887316971502630378744674713904894e-2\\ 
1.496263973303768103769644021861162037467544297982296706266438520924742416377446141864525561534462546 -0.8099445584437652165419405510818370160119788043007768997285904490810566734643868948837386818981491529e-3\\ 
1.594485336188642197940735516159224286083710846738877090026858573195349741382154659812520882993330767 -0.5746630658684697680649045635181315509631408189939541168152556274892970506064269313345328774429428320e-3\\ 
1.464644193717481857071842394969697625142729730498870114811186385942306964756565761319825750327407481 -0.5112551239641374366946387791388329788900864830045700177119732853602013884088666481660915156886471183e-3\\ 
1.442926823187799448886556676279581250768419006322026220052130747188952660542826665437394400252385619 -0.2440741964975975791511235750514365640785447542327275550201176592540885657391596515475714159553035966e-3\\
};
\addlegendentry{Model 4};

\end{axis}
\end{tikzpicture}%

%% file: case300-Model4-VC.tikz
%
%
\definecolor{mycolor1}{rgb}{0.00000,0.75000,0.75000}%
\definecolor{mycolor2}{rgb}{0.75000,0.00000,0.75000}%
\begin{tikzpicture}

\begin{axis}[%
width=5.705521cm,
height=6cm,
at={(0cm,0cm)},
scale only axis,
separate axis lines,
every outer x axis line/.append style={black},
every x tick label/.append style={font=\color{black}},
xmin=0,
xmax=16,
xlabel={n},
every outer y axis line/.append style={black},
every y tick label/.append style={font=\color{black}},
ymode=log,
ymin=0.01,
ymax=100000000000000,
yminorticks=true,
ylabel={$\text{|V}_\text{i}\text{[n]|}$},
legend style={at={(0.03,0.97)},anchor=north west,legend cell align=left,align=left,draw=black}
]
\addplot [color=blue,solid]
  table[row sep=crcr]{%
1	1.0507\\
2	1.78573337144854\\
3	8.54349424597091\\
4	59.2383195261587\\
5	414.62602118046\\
6	3061.34398928221\\
7	23310.1236911037\\
8	181832.984068686\\
9	1450115.63428621\\
10	11795112.7578521\\
11	97614542.141465\\
12	820045188.133284\\
13	6978597792.95293\\
14	60050743947.2215\\
15	521697128064.363\\
16	4569897553745.75\\
};
\addlegendentry{Maximum};

\addplot [color=black!50!green,solid]
  table[row sep=crcr]{%
1	1.0507\\
2	0.976172189604755\\
3	6.62203570604817\\
4	46.0594235750853\\
5	332.046427999532\\
6	2423.00600197388\\
7	18007.3658850335\\
8	139300.248656055\\
9	1089219.58858377\\
10	8800935.19272898\\
11	72934076.1942678\\
12	608310636.911081\\
13	5126863481.56396\\
14	43760416480.5669\\
15	377574741401.446\\
16	3291920714719.8\\
};
\addlegendentry{P75};

\addplot [color=red,solid]
  table[row sep=crcr]{%
1	1.0507\\
2	0.447380499765504\\
3	3.91133958389279\\
4	27.254341332429\\
5	196.944927354522\\
6	1462.05684476573\\
7	11106.1846319164\\
8	85826.1075904448\\
9	678872.745772251\\
10	5486252.42219415\\
11	45204214.4785529\\
12	378969216.96901\\
13	3233738494.92098\\
14	27786491122.2602\\
15	241094404987.336\\
16	2109886549410.99\\
};
\addlegendentry{Median};

\addplot [color=mycolor1,solid]
  table[row sep=crcr]{%
1	1.0507\\
2	0.278440420282762\\
3	1.16207120896263\\
4	8.00290618118128\\
5	55.9060495596933\\
6	408.576292386823\\
7	3158.80723726727\\
8	24675.753474845\\
9	193955.939452118\\
10	1562459.72844657\\
11	12848664.0829054\\
12	107476998.52918\\
13	911804238.993837\\
14	7826867504.85775\\
15	67852199681.5939\\
16	593189277854.899\\
};
\addlegendentry{P25};

\addplot [color=mycolor2,solid]
  table[row sep=crcr]{%
1	1.0507\\
2	0.0115082610636689\\
3	0.0332509289879613\\
4	0.247734222090887\\
5	1.4803237987959\\
6	1.58342009923456\\
7	37.6006449868811\\
8	540.268785883471\\
9	5487.33525191632\\
10	50121.2633748678\\
11	440689.012777718\\
12	3836407.21891851\\
13	33494293.1938391\\
14	294850418.887136\\
15	2620950702.51063\\
16	23515246025.1462\\
};
\addlegendentry{Minimum};

\end{axis}
\end{tikzpicture}%

%% file: case300-Model5-VC.tikz
%
%
\definecolor{mycolor1}{rgb}{0.00000,0.75000,0.75000}%
\definecolor{mycolor2}{rgb}{0.75000,0.00000,0.75000}%
\begin{tikzpicture}

\begin{axis}[%
width=5.705521cm,
height=6cm,
at={(0cm,0cm)},
scale only axis,
separate axis lines,
every outer x axis line/.append style={black},
every x tick label/.append style={font=\color{black}},
xmin=0,
xmax=16,
xlabel={n},
every outer y axis line/.append style={black},
every y tick label/.append style={font=\color{black}},
ymode=log,
ymin=1e-05,
ymax=10,
yminorticks=true,
ylabel={$\text{|V}_\text{i}\text{[n]|}$},
legend style={legend cell align=left,align=left,draw=black}
]
\addplot [color=blue,solid]
  table[row sep=crcr]{%
1	1.0507\\
2	0.962123670892302\\
3	0.570946933986721\\
4	0.324850574225532\\
5	0.126267667708403\\
6	0.0574402900299079\\
7	0.0198347483697719\\
8	0.0114478332303197\\
9	0.00577283624388799\\
10	0.00364793990014912\\
11	0.00218130324661038\\
12	0.00135806865712426\\
13	0.000911282673122756\\
14	0.000568310706202565\\
15	0.000404017746781264\\
16	0.000260311796328968\\
};
\addlegendentry{Maximum};

\addplot [color=black!50!green,solid]
  table[row sep=crcr]{%
1	1.0507\\
2	0.298339212901316\\
3	0.331884006882187\\
4	0.0934979543810484\\
5	0.0634562837276923\\
6	0.0198594836974029\\
7	0.0142702640767536\\
8	0.00633942768966751\\
9	0.00484776738178927\\
10	0.00239560575665715\\
11	0.00184269719393478\\
12	0.00103640240128104\\
13	0.000773889245972914\\
14	0.000467255553833046\\
15	0.00034300151127849\\
16	0.000219236724626264\\
};
\addlegendentry{P75};

\addplot [color=red,solid]
  table[row sep=crcr]{%
1	1.0507\\
2	0.159822735074156\\
3	0.287110958959988\\
4	0.0413785827691101\\
5	0.0528392931402765\\
6	0.0121810813130624\\
7	0.012424018354557\\
8	0.00435650224582919\\
9	0.00418510532570158\\
10	0.0017927911956878\\
11	0.00159726903839867\\
12	0.000788857267623788\\
13	0.000667770017283306\\
14	0.000360341936025131\\
15	0.000295411349132402\\
16	0.000170531345603455\\
};
\addlegendentry{Median};

\addplot [color=mycolor1,solid]
  table[row sep=crcr]{%
1	1.0507\\
2	0.0878116440336374\\
3	0.272390968051595\\
4	0.0235278988746628\\
5	0.047000187904551\\
6	0.00734107945659239\\
7	0.0104949452249105\\
8	0.00321514177605478\\
9	0.00364370716202216\\
10	0.00139597211808815\\
11	0.00141958162125028\\
12	0.000652618668971535\\
13	0.000604662126726866\\
14	0.000309581422633722\\
15	0.000269192473472735\\
16	0.000150552254538047\\
};
\addlegendentry{P25};

\addplot [color=mycolor2,solid]
  table[row sep=crcr]{%
1	1.0507\\
2	0.00421122954403965\\
3	0.0441262211562429\\
4	0.00230418616144165\\
5	0.00730674715039325\\
6	0.000822299928893117\\
7	0.0016912377063879\\
8	0.000416069590090952\\
9	0.000571653075292509\\
10	0.000190644288894742\\
11	0.00021567430248223\\
12	9.02865866207725e-05\\
13	8.94506614179379e-05\\
14	4.33132057686024e-05\\
15	3.9392940183456e-05\\
16	2.11546324907903e-05\\
};
\addlegendentry{Minimum};

\end{axis}
\end{tikzpicture}%

%% file: case39PQ-HELMBusVoltages.tikz
%
%
\definecolor{mycolor1}{rgb}{0.00000,0.75000,0.75000}%
\definecolor{mycolor2}{rgb}{0.75000,0.00000,0.75000}%
\begin{tikzpicture}

\begin{axis}[%
width=5.705521cm,
height=6cm,
at={(0cm,0cm)},
scale only axis,
separate axis lines,
every outer x axis line/.append style={black},
every x tick label/.append style={font=\color{black}},
xmin=0,
xmax=1,
xlabel={z Value},
every outer y axis line/.append style={black},
every y tick label/.append style={font=\color{black}},
ymin=0,
ymax=2.5,
ylabel={Voltage Magnitude},
title={case39PQ-Bus Voltages from HELM},
legend style={at={(0.97,0.03)},anchor=south east,legend cell align=left,align=left,draw=black}
]
\addplot [color=blue,solid]
  table[row sep=crcr]{%
0	2.13014715911593\\
0.01	2.12890931384455\\
0.02	2.12764079701938\\
0.03	2.1263414363258\\
0.04	2.12501105144181\\
0.05	2.12364945386693\\
0.06	2.12225644674175\\
0.07	2.12083182465757\\
0.08	2.11937537345569\\
0.09	2.11788687001588\\
0.1	2.11636608203356\\
0.11	2.11481276778503\\
0.12	2.11322667588023\\
0.13	2.11160754500239\\
0.14	2.10995510363392\\
0.15	2.10826906976777\\
0.16	2.10654915060344\\
0.17	2.10479504222704\\
0.18	2.10300642927426\\
0.19	2.10118298457543\\
0.2	2.09932436878167\\
0.21	2.09743022997101\\
0.22	2.09550020323338\\
0.23	2.09353391023308\\
0.24	2.09153095874761\\
0.25	2.08949094218119\\
0.26	2.08741343905166\\
0.27	2.08529801244884\\
0.28	2.08314420946281\\
0.29	2.08095156058\\
0.3	2.07871957904507\\
0.31	2.07644776018638\\
0.32	2.07413558070247\\
0.33	2.07178249790715\\
0.34	2.0693879489301\\
0.35	2.06695134987013\\
0.36	2.06447209489755\\
0.37	2.06194955530222\\
0.38	2.05938307848315\\
0.39	2.0567719868755\\
0.4	2.05411557681021\\
0.41	2.05141311730126\\
0.42	2.04866384875489\\
0.43	2.04586698159479\\
0.44	2.04302169479655\\
0.45	2.0401271343242\\
0.46	2.03718241146063\\
0.47	2.03418660102339\\
0.48	2.03113873945602\\
0.49	2.02803782278432\\
0.5	2.02488280442597\\
0.51	2.0216725928403\\
0.52	2.01840604900422\\
0.53	2.01508198369807\\
0.54	2.01169915458413\\
0.55	2.0082562630579\\
0.56	2.00475195085065\\
0.57	2.00118479635859\\
0.58	1.99755331067175\\
0.59	1.99385593327192\\
0.6	1.99009102736559\\
0.61	1.98625687481342\\
0.62	1.98235167061288\\
0.63	1.97837351688528\\
0.64	1.97432041631163\\
0.65	1.97019026495468\\
0.66	1.96598084439547\\
0.67	1.96168981310298\\
0.68	1.95731469694354\\
0.69	1.95285287872299\\
0.7	1.94830158663856\\
0.71	1.94365788149847\\
0.72	1.93891864254478\\
0.73	1.93408055168842\\
0.74	1.92914007593404\\
0.75	1.92409344773379\\
0.76	1.91893664296431\\
0.77	1.91366535616591\\
0.78	1.90827497261642\\
0.79	1.90276053673109\\
0.8	1.89711671618045\\
0.81	1.89133776099573\\
0.82	1.88541745678209\\
0.83	1.87934907098115\\
0.84	1.87312529094672\\
0.85	1.86673815263161\\
0.86	1.86017896029011\\
0.87	1.85343821284167\\
0.88	1.84650589004622\\
0.89	1.83936860198885\\
0.9	1.8320155989253\\
0.91	1.82443184163827\\
0.92	1.81660137402086\\
0.93	1.80850626149182\\
0.94	1.80012622746159\\
0.95	1.79143820571363\\
0.96	1.78241577433401\\
0.97	1.77302843186477\\
0.98	1.76324065907384\\
0.99	1.75301068353058\\
1	1.74228882327307\\
};
\addlegendentry{Maximum};

\addplot [color=black!50!green,solid]
  table[row sep=crcr]{%
0	1.96081535825113\\
0.01	1.96008249424409\\
0.02	1.95932157648119\\
0.03	1.95853246766243\\
0.04	1.95771502325926\\
0.05	1.95686909137156\\
0.06	1.95599451257611\\
0.07	1.95509111976644\\
0.08	1.95383987822623\\
0.09	1.95250889945521\\
0.1	1.95119698520909\\
0.11	1.95011136386899\\
0.12	1.9489960888814\\
0.13	1.9478509409171\\
0.14	1.94667569158119\\
0.15	1.94547010316511\\
0.16	1.94423392838525\\
0.17	1.94296691010734\\
0.18	1.9416151073287\\
0.19	1.94012673266592\\
0.2	1.93860749706275\\
0.21	1.93705710035656\\
0.22	1.93547523113663\\
0.23	1.93386156636984\\
0.24	1.93221577100584\\
0.25	1.93053749756057\\
0.26	1.92882638567677\\
0.27	1.92708206165996\\
0.28	1.92530413798853\\
0.29	1.92349221279614\\
0.3	1.92164586932476\\
0.31	1.91976467534627\\
0.32	1.9178481825507\\
0.33	1.9158883304958\\
0.34	1.9138880336041\\
0.35	1.91185092375884\\
0.36	1.90977648010723\\
0.37	1.90766416202488\\
0.38	1.90551340824697\\
0.39	1.9033236359454\\
0.4	1.90109423974807\\
0.41	1.89882459069581\\
0.42	1.89651403513228\\
0.43	1.89416189352148\\
0.44	1.89176745918736\\
0.45	1.88932999696905\\
0.46	1.88684874178501\\
0.47	1.88432289709855\\
0.48	1.88175163327629\\
0.49	1.87913408583066\\
0.5	1.87646935353621\\
0.51	1.8737564964086\\
0.52	1.87099453353412\\
0.53	1.86818244073601\\
0.54	1.86531914806239\\
0.55	1.86240353707919\\
0.56	1.85943443794918\\
0.57	1.85641062627626\\
0.58	1.85333081969174\\
0.59	1.8501936741564\\
0.6	1.84699777994898\\
0.61	1.8437416573082\\
0.62	1.84054891566589\\
0.63	1.83732684632951\\
0.64	1.83403911763645\\
0.65	1.83068391224169\\
0.66	1.82725931373873\\
0.67	1.82376329920324\\
0.68	1.82019373099362\\
0.69	1.81654834771642\\
0.7	1.81282475425031\\
0.71	1.80902041070601\\
0.72	1.80513262018019\\
0.73	1.80115851513873\\
0.74	1.79709504223697\\
0.75	1.79293894535229\\
0.76	1.7886867465649\\
0.77	1.78433472477552\\
0.78	1.77987889159094\\
0.79	1.77531496403866\\
0.8	1.7706383335857\\
0.81	1.76584403083054\\
0.82	1.76092668510371\\
0.83	1.75588047803412\\
0.84	1.75069908980139\\
0.85	1.74537563459853\\
0.86	1.73990260553838\\
0.87	1.73427174539623\\
0.88	1.72847401730846\\
0.89	1.72249940987139\\
0.9	1.71633684453451\\
0.91	1.70997399232895\\
0.92	1.7033970684874\\
0.93	1.69659058047634\\
0.94	1.68953701692284\\
0.95	1.68221645997732\\
0.96	1.67460609647995\\
0.97	1.666679593531\\
0.98	1.65840628942494\\
0.99	1.64975012847939\\
1	1.6406682330909\\
};
\addlegendentry{P75};

\addplot [color=red,solid]
  table[row sep=crcr]{%
0	1.92448668656981\\
0.01	1.92325995293063\\
0.02	1.92200835346594\\
0.03	1.92073174265509\\
0.04	1.91942996851172\\
0.05	1.91810287244475\\
0.06	1.91675028911147\\
0.07	1.91537204626279\\
0.08	1.91396796457995\\
0.09	1.91253785750252\\
0.1	1.91108153104716\\
0.11	1.90959878361676\\
0.12	1.9080894057994\\
0.13	1.90655318015675\\
0.14	1.90498988100125\\
0.15	1.90339927416156\\
0.16	1.90178111673563\\
0.17	1.90013515683073\\
0.18	1.89846113328979\\
0.19	1.89685191883805\\
0.2	1.89528688716238\\
0.21	1.89369212689519\\
0.22	1.89206733785545\\
0.23	1.89041220876799\\
0.24	1.88883962735343\\
0.25	1.88736060973007\\
0.26	1.88584978473293\\
0.27	1.88430679653471\\
0.28	1.88273127600059\\
0.29	1.88112284019162\\
0.3	1.87948109183985\\
0.31	1.87780561879343\\
0.32	1.87609599342957\\
0.33	1.87435177203329\\
0.34	1.87257249413954\\
0.35	1.87075768183618\\
0.36	1.86890683902514\\
0.37	1.86701945063856\\
0.38	1.865094981807\\
0.39	1.86313287697577\\
0.4	1.86113255896596\\
0.41	1.85909342797558\\
0.42	1.85701486051658\\
0.43	1.85489620828247\\
0.44	1.85273679694126\\
0.45	1.85053592484758\\
0.46	1.84829286166753\\
0.47	1.84600684690893\\
0.48	1.84367708834915\\
0.49	1.8413027603516\\
0.5	1.83888300206141\\
0.51	1.83641691546955\\
0.52	1.83390356333363\\
0.53	1.83134196694229\\
0.54	1.82873110370874\\
0.55	1.82606990457725\\
0.56	1.82335725122487\\
0.57	1.82059197303796\\
0.58	1.81777284384164\\
0.59	1.81481727556474\\
0.6	1.81165673052124\\
0.61	1.80843837535915\\
0.62	1.80516071429055\\
0.63	1.80182217473708\\
0.64	1.79842110191294\\
0.65	1.79495575290183\\
0.66	1.79142429016899\\
0.67	1.78782477444113\\
0.68	1.78415515687741\\
0.69	1.78041327044355\\
0.7	1.77659682038746\\
0.71	1.77270337369976\\
0.72	1.76873034742361\\
0.73	1.76467499565684\\
0.74	1.76053439506276\\
0.75	1.7563054286759\\
0.76	1.7519847677501\\
0.77	1.74756885135252\\
0.78	1.74305386335159\\
0.79	1.73843570638041\\
0.8	1.73370997227519\\
0.81	1.72887190838747\\
0.82	1.72391637904357\\
0.83	1.71883782126729\\
0.84	1.7136301936743\\
0.85	1.70828691708861\\
0.86	1.70280080407685\\
0.87	1.69716398549423\\
0.88	1.6913678204862\\
0.89	1.68540272057473\\
0.9	1.67925812166161\\
0.91	1.67292225614249\\
0.92	1.66638196052238\\
0.93	1.65962244095821\\
0.94	1.65256293628929\\
0.95	1.64519170132436\\
0.96	1.63754557733729\\
0.97	1.62959865392945\\
0.98	1.62132186568239\\
0.99	1.61268104136802\\
1	1.60363551550618\\
};
\addlegendentry{Median};

\addplot [color=mycolor1,solid]
  table[row sep=crcr]{%
0	1.76467218890732\\
0.01	1.76363232349267\\
0.02	1.76257092279995\\
0.03	1.76148786198343\\
0.04	1.76038301065062\\
0.05	1.75925623274174\\
0.06	1.75810738640273\\
0.07	1.75693632385136\\
0.08	1.75574289123625\\
0.09	1.75452692848837\\
0.1	1.75328826916468\\
0.11	1.75202674028356\\
0.12	1.75074216215159\\
0.13	1.74943434818126\\
0.14	1.74810310469911\\
0.15	1.74674823074387\\
0.16	1.74536951785401\\
0.17	1.74396674984409\\
0.18	1.74253970256946\\
0.19	1.74108814367849\\
0.2	1.73961183235168\\
0.21	1.73811051902698\\
0.22	1.73658394511038\\
0.23	1.73503184267107\\
0.24	1.73345393412007\\
0.25	1.7318499318715\\
0.26	1.73021953798527\\
0.27	1.72856244379019\\
0.28	1.72687832948616\\
0.29	1.72516686372408\\
0.3	1.72342770316219\\
0.31	1.72166049199711\\
0.32	1.71986486146799\\
0.33	1.71804042933199\\
0.34	1.71618679930906\\
0.35	1.71430356049382\\
0.36	1.71239028673249\\
0.37	1.71044653596202\\
0.38	1.70847184950901\\
0.39	1.70646575134521\\
0.4	1.70442774729657\\
0.41	1.70235732420215\\
0.42	1.70025394901924\\
0.43	1.69811706787028\\
0.44	1.69594610502722\\
0.45	1.6937404618281\\
0.46	1.69149951552049\\
0.47	1.68922261802562\\
0.48	1.68690909461664\\
0.49	1.68455824250363\\
0.5	1.6821693293173\\
0.51	1.67974159148243\\
0.52	1.67727423247126\\
0.53	1.67476642092576\\
0.54	1.67221728863676\\
0.55	1.66962592836646\\
0.56	1.66699139149913\\
0.57	1.66431268550354\\
0.58	1.66158877118806\\
0.59	1.65881855972772\\
0.6	1.65600090943964\\
0.61	1.65313462228025\\
0.62	1.6502184400347\\
0.63	1.6472510401646\\
0.64	1.64423103127623\\
0.65	1.6411569481657\\
0.66	1.63802724639214\\
0.67	1.63484029632273\\
0.68	1.63159437658543\\
0.69	1.62828766685584\\
0.7	1.62491823989373\\
0.71	1.62148405273144\\
0.72	1.61798293690125\\
0.73	1.61441258757055\\
0.74	1.61077055143171\\
0.75	1.60705421316788\\
0.76	1.60326078028448\\
0.77	1.59938726605846\\
0.78	1.59543047031206\\
0.79	1.59138695766162\\
0.8	1.58725303282405\\
0.81	1.58302471247942\\
0.82	1.57869769308346\\
0.83	1.57426731389365\\
0.84	1.56972851430599\\
0.85	1.56507578436354\\
0.86	1.56030310576749\\
0.87	1.55540389497444\\
0.88	1.55037089765312\\
0.89	1.54519611764719\\
0.9	1.53987070298474\\
0.91	1.53438481341307\\
0.92	1.52872743442314\\
0.93	1.52288619411318\\
0.94	1.5168471147364\\
0.95	1.51059430327366\\
0.96	1.5041095625628\\
0.97	1.49737189493716\\
0.98	1.49035685928237\\
0.99	1.48303572476253\\
1	1.47537433658545\\
};
\addlegendentry{P25};

\addplot [color=mycolor2,solid]
  table[row sep=crcr]{%
0	1.58434793901313\\
0.01	1.58379650360711\\
0.02	1.58322634256086\\
0.03	1.58263735931989\\
0.04	1.58202945261012\\
0.05	1.58140251633874\\
0.06	1.58075643948968\\
0.07	1.58009110601336\\
0.08	1.57940639471059\\
0.09	1.57870217911024\\
0.1	1.5779783273404\\
0.11	1.5772347019928\\
0.12	1.57647115997992\\
0.13	1.57568755238469\\
0.14	1.5748837243022\\
0.15	1.5740595146731\\
0.16	1.57321475610819\\
0.17	1.57234927470373\\
0.18	1.57146288984703\\
0.19	1.57055541401169\\
0.2	1.56962665254193\\
0.21	1.56867640342536\\
0.22	1.56770445705368\\
0.23	1.56671059597027\\
0.24	1.56569459460429\\
0.25	1.56465621899014\\
0.26	1.5635952264715\\
0.27	1.56251136538911\\
0.28	1.56140437475109\\
0.29	1.56027398388471\\
0.3	1.55911991206852\\
0.31	1.55794186814347\\
0.32	1.55673955010158\\
0.33	1.55551264465076\\
0.34	1.55426082675404\\
0.35	1.55298375914148\\
0.36	1.55168109179291\\
0.37	1.5503524613892\\
0.38	1.54899749073013\\
0.39	1.54761578811599\\
0.4	1.54620694669061\\
0.41	1.54477054374264\\
0.42	1.54330613996201\\
0.43	1.54181327864799\\
0.44	1.54029148486512\\
0.45	1.53874026454285\\
0.46	1.53715910351413\\
0.47	1.5355474664881\\
0.48	1.53390479595139\\
0.49	1.53223051099156\\
0.5	1.53052400603645\\
0.51	1.52878464950166\\
0.52	1.52701178233796\\
0.53	1.5252047164692\\
0.54	1.52336273311108\\
0.55	1.52148508095664\\
0.56	1.51957097421209\\
0.57	1.51761959043479\\
0.58	1.51563006765798\\
0.59	1.51360150932435\\
0.6	1.51153295839658\\
0.61	1.50942342790388\\
0.62	1.50727187325878\\
0.63	1.50507719719601\\
0.64	1.50283824520726\\
0.65	1.50055380135737\\
0.66	1.49822258375811\\
0.67	1.49584323961181\\
0.68	1.49341433976401\\
0.69	1.49093437270345\\
0.7	1.48840173793807\\
0.71	1.48581473866574\\
0.72	1.48317157364579\\
0.73	1.48047032816204\\
0.74	1.47770896395039\\
0.75	1.47488530794164\\
0.76	1.47199703964546\\
0.77	1.4690416769688\\
0.78	1.46601656022494\\
0.79	1.46291883404291\\
0.8	1.45974542682987\\
0.81	1.45649302737006\\
0.82	1.45315805805805\\
0.83	1.44973664416445\\
0.84	1.44622457844832\\
0.85	1.44261728060685\\
0.86	1.43890974243915\\
0.87	1.43509644982594\\
0.88	1.43117145909518\\
0.89	1.4271280736142\\
0.9	1.42295893683844\\
0.91	1.4186558580382\\
0.92	1.41420968126798\\
0.93	1.4096101193637\\
0.94	1.40484554833608\\
0.95	1.39990275106055\\
0.96	1.39476659409794\\
0.97	1.38941961493467\\
0.98	1.38384148725429\\
0.99	1.37800831704174\\
1	1.3718916990735\\
};
\addlegendentry{Minimum};

\end{axis}
\end{tikzpicture}%

%% file: case39PQ-MatPowerBusVoltages.tikz
%
%
\definecolor{mycolor1}{rgb}{0.00000,0.75000,0.75000}%
\definecolor{mycolor2}{rgb}{0.75000,0.00000,0.75000}%
\begin{tikzpicture}

\begin{axis}[%
width=5.705521cm,
height=6cm,
at={(0cm,0cm)},
scale only axis,
separate axis lines,
every outer x axis line/.append style={black},
every x tick label/.append style={font=\color{black}},
xmin=0,
xmax=1,
xlabel={Load Factor},
every outer y axis line/.append style={black},
every y tick label/.append style={font=\color{black}},
ymin=0,
ymax=2.5,
ylabel={Voltage Magnitude},
title={case39PQ-Bus Voltages from MATPOWER},
legend style={at={(0.97,0.6)},anchor=south east,legend cell align=left,align=left,draw=black}
]
\addplot [color=blue,solid]
  table[row sep=crcr]{%
1	1.06359999970567\\
0.99	1.0491731431323\\
0.98	1.03517050366973\\
0.97	1.02154984745351\\
0.96	1.00827588025041\\
0.95	0.995318824523899\\
0.94	0.982653362211996\\
0.93	0.970706148099589\\
0.92	0.962097025127069\\
0.91	0.953698095856898\\
0.9	0.945497397950709\\
0.89	0.937484486850807\\
0.88	0.929650239345802\\
0.87	0.921986697891935\\
0.86	0.914486948123954\\
0.85	0.907504434338287\\
0.84	0.902007158393027\\
0.83	0.896627036798533\\
0.82	0.891361574228249\\
0.81	0.886208914566227\\
0.8	0.881167855262175\\
0.79	0.876237880428796\\
0.78	0.871419215590727\\
0.77	0.866712908414408\\
0.76	0.8621209415913\\
0.75	0.857646386499184\\
0.74	0.853293609578474\\
0.73	0.849068547865452\\
0.72	0.844979076244842\\
0.71	0.841035497202057\\
0.7	0.837251194576157\\
0.69	0.833643505951335\\
0.68	0.83023488221097\\
0.67	0.827054411516553\\
0.66	0.824139772379179\\
0.65	0.821539607319989\\
0.64	0.81931608929149\\
0.63	0.817546926855606\\
0.62	0.816325002641064\\
0.61	0.815752242920025\\
0.6	0.815923216538148\\
0.59	0.816896674575229\\
0.58	0.81866348283506\\
0.57	0.82113167325901\\
0.56	0.824144862460213\\
0.55	0.827524524181843\\
0.54	0.831109176568954\\
0.53	0.834773808267114\\
0.52	0.838431459478922\\
0.51	0.842026169890497\\
0.5	0.845524266765048\\
0.49	0.848907022492324\\
0.48	0.852165321619063\\
0.47	0.855296056335687\\
0.46	0.858299778928168\\
0.45	0.861179202839337\\
0.44	0.863938256638141\\
0.43	0.866581492303444\\
0.42	0.869113719010444\\
0.41	0.871539780154659\\
0.4	0.873864421362358\\
0.39	0.876092216311566\\
0.38	0.878227529256172\\
0.37	0.880274500789742\\
0.36	0.882237048248517\\
0.35	0.88411887526355\\
0.34	0.885923486976542\\
0.33	0.88765420873476\\
0.32	0.889314206931648\\
0.31	0.890906511222398\\
0.3	0.892434037719219\\
0.29	0.893899613025769\\
0.28	0.895305999145252\\
0.27	0.896655919419107\\
0.26	0.897952085739014\\
0.25	0.899197227333004\\
0.24	0.900394121459624\\
0.23	0.9015456263502\\
0.22	0.90265471670992\\
0.21	0.903724522010313\\
0.2	0.904758367655297\\
0.19	0.905759818850962\\
0.18	0.906732726611166\\
0.17	0.907681274735841\\
0.16	0.908610025742909\\
0.15	0.909523962557731\\
0.14	0.910428521222629\\
0.13	0.91132960799548\\
0.12	0.912233592084884\\
0.11	0.913147263235065\\
0.1	0.914077742021654\\
0.09	0.915032330977251\\
0.08	0.916018297730893\\
0.07	0.917042588393478\\
0.06	0.918111480999961\\
0.05	0.919230204007613\\
0.04	0.920402560561692\\
0.03	0.921630610340993\\
0.02	0.922914461648095\\
0.01	0.924252213634804\\
0	0.925640063641119\\
};
\addlegendentry{Maximum};

\addplot [color=black!50!green,solid]
  table[row sep=crcr]{%
1	1.0489470561573\\
0.99	1.03361276793891\\
0.98	1.01763447874303\\
0.97	1.0054916065198\\
0.96	0.993236645696863\\
0.95	0.980943396858377\\
0.94	0.9689460117886\\
0.93	0.958113222329903\\
0.92	0.948449418204242\\
0.91	0.938662717452465\\
0.9	0.930374035737549\\
0.89	0.921506681633004\\
0.88	0.911634049616256\\
0.87	0.903022421892074\\
0.86	0.894427032488905\\
0.85	0.885742571015732\\
0.84	0.878137480199315\\
0.83	0.871166042576856\\
0.82	0.864340382508339\\
0.81	0.857657974740122\\
0.8	0.851117109428739\\
0.79	0.844716929983089\\
0.78	0.838652888817813\\
0.77	0.833298610307178\\
0.76	0.82806891175507\\
0.75	0.822967351046772\\
0.74	0.817998931909279\\
0.73	0.813170411689796\\
0.72	0.808490710745865\\
0.71	0.803971457305287\\
0.7	0.799627713297317\\
0.69	0.795478940796366\\
0.68	0.79155028337848\\
0.67	0.787874245124941\\
0.66	0.784492833977564\\
0.65	0.781460151835321\\
0.64	0.778845166349454\\
0.63	0.776733812288343\\
0.62	0.775228407637159\\
0.61	0.774440624609228\\
0.6	0.774473106501018\\
0.59	0.775387937066051\\
0.58	0.777171550917476\\
0.57	0.779719032988761\\
0.56	0.782855360935256\\
0.55	0.786382522139917\\
0.54	0.790012057913574\\
0.53	0.793706442501859\\
0.52	0.797383692220647\\
0.51	0.800983096669666\\
0.5	0.804468526019236\\
0.49	0.807820241598006\\
0.48	0.811028982899759\\
0.47	0.814091990910144\\
0.46	0.817101934022207\\
0.45	0.820305938679731\\
0.44	0.823363500936624\\
0.43	0.826280384850114\\
0.42	0.829062614155618\\
0.41	0.831999127865137\\
0.4	0.834866507813003\\
0.39	0.837616685694953\\
0.38	0.840254812119682\\
0.37	0.842785808306609\\
0.36	0.845214358461314\\
0.35	0.84754491533744\\
0.34	0.849781714253173\\
0.33	0.851928792549274\\
0.32	0.8539900126095\\
0.31	0.855969087320357\\
0.3	0.857869607356573\\
0.29	0.859695070028005\\
0.28	0.86144890966636\\
0.27	0.863134529701802\\
0.26	0.864755336701412\\
0.25	0.866314776725219\\
0.24	0.867816374405294\\
0.23	0.869263775166175\\
0.22	0.870660790970098\\
0.21	0.872011449870181\\
0.2	0.873320049458463\\
0.19	0.874591213964313\\
0.18	0.875829954234094\\
0.17	0.877041729037424\\
0.16	0.878232505015121\\
0.15	0.879408811029682\\
0.14	0.880577780643803\\
0.13	0.881747173952876\\
0.12	0.882925367198063\\
0.11	0.884121295904706\\
0.1	0.885344335513202\\
0.09	0.886604103821818\\
0.08	0.887910173625006\\
0.07	0.889271693257403\\
0.06	0.890696928063781\\
0.05	0.892192755887101\\
0.04	0.893764170418408\\
0.03	0.895413860898066\\
0.02	0.897141937750926\\
0.01	0.898678785943715\\
0	0.899997403035611\\
};
\addlegendentry{P75};

\addplot [color=red,solid]
  table[row sep=crcr]{%
1	1.03114015357927\\
0.99	1.01828573073889\\
0.98	1.00561427348143\\
0.97	0.99298243760364\\
0.96	0.98146188147822\\
0.95	0.969901751787889\\
0.94	0.958286591724827\\
0.93	0.947899402347254\\
0.92	0.937674826813653\\
0.91	0.928777052150261\\
0.9	0.919337950498288\\
0.89	0.91068742584133\\
0.88	0.90035611189201\\
0.87	0.891226741219085\\
0.86	0.88180914918574\\
0.85	0.872164117433182\\
0.84	0.862633297781601\\
0.83	0.853261786066325\\
0.82	0.844046186798502\\
0.81	0.834939734825345\\
0.8	0.825245659619694\\
0.79	0.815717445818493\\
0.78	0.806358036204465\\
0.77	0.797172085235196\\
0.76	0.788166217426362\\
0.75	0.779349367606885\\
0.74	0.770733226182688\\
0.73	0.7623328206114\\
0.72	0.754167274934335\\
0.71	0.746260802949588\\
0.7	0.738704680702739\\
0.69	0.731741455398732\\
0.68	0.725247039125217\\
0.67	0.719162987198881\\
0.66	0.713562022511051\\
0.65	0.708536163309184\\
0.64	0.703580488411379\\
0.63	0.698824061057594\\
0.62	0.695738935215365\\
0.61	0.694339634305864\\
0.6	0.694336033627245\\
0.59	0.695469501577606\\
0.58	0.696501442197259\\
0.57	0.700121497292813\\
0.56	0.704211512753973\\
0.55	0.708559898858351\\
0.54	0.714120922235828\\
0.53	0.71969412602542\\
0.52	0.724853851942927\\
0.51	0.729350889615527\\
0.5	0.733573751357245\\
0.49	0.738367717859262\\
0.48	0.742884769341783\\
0.47	0.747129559010693\\
0.46	0.751111371582249\\
0.45	0.754841773767245\\
0.44	0.758333218338265\\
0.43	0.761598234339948\\
0.42	0.764648974271855\\
0.41	0.767496977155681\\
0.4	0.770081245231788\\
0.39	0.772409901159803\\
0.38	0.774572386959085\\
0.37	0.776576869539456\\
0.36	0.778430918448002\\
0.35	0.780141557054007\\
0.34	0.781715316596272\\
0.33	0.783158290730965\\
0.32	0.784476189374827\\
0.31	0.785674391366952\\
0.3	0.786757995933128\\
0.29	0.7877318732309\\
0.28	0.788600714443755\\
0.27	0.789369082017226\\
0.26	0.790041460711882\\
0.25	0.790844368540324\\
0.24	0.79162379224266\\
0.23	0.792323232413507\\
0.22	0.792947502885715\\
0.21	0.793501641544562\\
0.2	0.793990980842687\\
0.19	0.79442122273189\\
0.18	0.794798517227236\\
0.17	0.795129542855626\\
0.16	0.795421585852425\\
0.15	0.795682613052766\\
0.14	0.795921330901612\\
0.13	0.796147219889626\\
0.12	0.796370530199522\\
0.11	0.796602220920074\\
0.1	0.796853822808147\\
0.09	0.797137204760504\\
0.08	0.797464228861271\\
0.07	0.797846290057986\\
0.06	0.798293755166482\\
0.05	0.798815340626886\\
0.04	0.799417494314848\\
0.03	0.80010386543982\\
0.02	0.800914677834218\\
0.01	0.80218981817768\\
0	0.80355203766745\\
};
\addlegendentry{Median};

\addplot [color=mycolor1,solid]
  table[row sep=crcr]{%
1	1.01230948073916\\
0.99	1.00092332677381\\
0.98	0.991373101701824\\
0.97	0.980187017445553\\
0.96	0.969863765556398\\
0.95	0.960704743505011\\
0.94	0.951424582713799\\
0.93	0.941532803073559\\
0.92	0.930372007877063\\
0.91	0.919597784255069\\
0.9	0.908111240518014\\
0.89	0.898404368390187\\
0.88	0.888188962897661\\
0.87	0.87758758145969\\
0.86	0.867447310621986\\
0.85	0.857165003076205\\
0.84	0.847030086337682\\
0.83	0.837045092643744\\
0.82	0.827479614254387\\
0.81	0.818497516874711\\
0.8	0.80966040497771\\
0.79	0.800377641898785\\
0.78	0.791182105128447\\
0.77	0.782140781902018\\
0.76	0.773257207386263\\
0.75	0.764536702440859\\
0.74	0.755986712277951\\
0.73	0.747617256527074\\
0.72	0.739441526734912\\
0.71	0.731476680299172\\
0.7	0.723744896671767\\
0.69	0.71627478214329\\
0.68	0.709345423718964\\
0.67	0.702991719192651\\
0.66	0.697153686062818\\
0.65	0.691929751314135\\
0.64	0.687442707223516\\
0.63	0.683841821677355\\
0.62	0.680751126708679\\
0.61	0.678551535160121\\
0.6	0.677574863954277\\
0.59	0.67790828786587\\
0.58	0.679119645731483\\
0.57	0.679654019642219\\
0.56	0.681669718290718\\
0.55	0.685438367342518\\
0.54	0.689495309545287\\
0.53	0.693651324391567\\
0.52	0.697777399673335\\
0.51	0.701792453144166\\
0.5	0.705649207341627\\
0.49	0.709322636276793\\
0.48	0.712801744249766\\
0.47	0.716084104511548\\
0.46	0.719172360412908\\
0.45	0.722065275368406\\
0.44	0.724380620958342\\
0.43	0.726606402659741\\
0.42	0.728685528790697\\
0.41	0.730624204981851\\
0.4	0.732428538735741\\
0.39	0.734104431686914\\
0.38	0.735657524092948\\
0.37	0.737093172511386\\
0.36	0.738416448498989\\
0.35	0.739692335350773\\
0.34	0.740929733120008\\
0.33	0.742065831930574\\
0.32	0.743104906077189\\
0.31	0.74414261592299\\
0.3	0.745172387265272\\
0.29	0.746097856289045\\
0.28	0.746923459461147\\
0.27	0.747653511473682\\
0.26	0.748292253457321\\
0.25	0.748869590772039\\
0.24	0.749482402254586\\
0.23	0.750028030033989\\
0.22	0.7505109582499\\
0.21	0.750935926921611\\
0.2	0.75130799952706\\
0.19	0.751632635222084\\
0.18	0.751915764865112\\
0.17	0.75216386905556\\
0.16	0.752384055012945\\
0.15	0.752564098478586\\
0.14	0.75253841024781\\
0.13	0.752357598911135\\
0.12	0.752133776800708\\
0.11	0.751911599896131\\
0.1	0.751701795032075\\
0.09	0.751515423726952\\
0.08	0.751363563755022\\
0.07	0.751256882234927\\
0.06	0.751205113390803\\
0.05	0.751346098681125\\
0.04	0.751857069998305\\
0.03	0.752464215229462\\
0.02	0.753168053870823\\
0.01	0.753965887582581\\
0	0.754851948151287\\
};
\addlegendentry{P25};

\addplot [color=mycolor2,solid]
  table[row sep=crcr]{%
1	0.984099999860666\\
0.99	0.975259438099506\\
0.98	0.964754113932001\\
0.97	0.95218336639418\\
0.96	0.939938149556363\\
0.95	0.927990638517374\\
0.94	0.916317163314829\\
0.93	0.904897462719014\\
0.92	0.893714112984702\\
0.91	0.882752085975261\\
0.9	0.871998404904488\\
0.89	0.861441875250427\\
0.88	0.849802231614556\\
0.87	0.836490405923003\\
0.86	0.823309084956141\\
0.85	0.810244694437789\\
0.84	0.797284232097154\\
0.83	0.784415103935775\\
0.82	0.77162496890942\\
0.81	0.75890158722534\\
0.8	0.746232667715271\\
0.79	0.733605709788965\\
0.78	0.721007835333152\\
0.77	0.708425605633388\\
0.76	0.695844818008132\\
0.75	0.683250276454918\\
0.74	0.670625530401092\\
0.73	0.657952575999971\\
0.72	0.645211516038736\\
0.71	0.632380178824118\\
0.7	0.619433706043945\\
0.69	0.606344139582667\\
0.68	0.593080076898799\\
0.67	0.579606540645027\\
0.66	0.565885349477085\\
0.65	0.55187652741535\\
0.64	0.537541698187977\\
0.63	0.522850972783991\\
0.62	0.507795287308959\\
0.61	0.492405501454188\\
0.6	0.476775718512893\\
0.59	0.461079574341197\\
0.58	0.445559898429255\\
0.57	0.430480454473624\\
0.56	0.416059349445817\\
0.55	0.402424444029295\\
0.54	0.38961092499868\\
0.53	0.377587242094727\\
0.52	0.366285981172666\\
0.51	0.355626853799383\\
0.5	0.345529979379226\\
0.49	0.335922110113209\\
0.48	0.326738764741237\\
0.47	0.317924318568437\\
0.46	0.309431216908596\\
0.45	0.301218906215617\\
0.44	0.29325275777628\\
0.43	0.285503095517858\\
0.42	0.277944361321362\\
0.41	0.270554416552826\\
0.4	0.263313965027564\\
0.39	0.256206078952305\\
0.38	0.249215809965613\\
0.37	0.242329869559668\\
0.36	0.235536365704801\\
0.35	0.228824584901668\\
0.34	0.222184810971892\\
0.33	0.215608173633514\\
0.32	0.209086521320267\\
0.31	0.202612313842118\\
0.3	0.19617853140023\\
0.29	0.189778597209429\\
0.28	0.183406311585951\\
0.27	0.177055795861993\\
0.26	0.170721444919727\\
0.25	0.164397887520571\\
0.24	0.15807995396144\\
0.23	0.151762650936892\\
0.22	0.145441143840984\\
0.21	0.139110747118823\\
0.2	0.13276692368559\\
0.19	0.126405294873573\\
0.18	0.120021662838191\\
0.17	0.113612047825097\\
0.16	0.107172743115119\\
0.15	0.100700390717786\\
0.14	0.0941920808081343\\
0.13	0.0876454772432429\\
0.12	0.0810589699150628\\
0.11	0.074431851797436\\
0.1	0.0677645139664483\\
0.09	0.061058645482854\\
0.08	0.0543174172264386\\
0.07	0.0475456208613943\\
0.06	0.0407497285039544\\
0.05	0.0339378386584969\\
0.04	0.0271194827559918\\
0.03	0.0203052854318251\\
0.02	0.0135064980589912\\
0.01	0.00673445233716585\\
0	1.49194347803321e-14\\
};
\addlegendentry{Minimum};

\end{axis}
\end{tikzpicture}%